\documentclass[a4paper,11pt,leqno]{article}
\usepackage{amsmath, amssymb, amsthm}
\usepackage{graphicx} 

\numberwithin{equation}{section}

\begin{document}

\def\S{{\mathbb S}}\def\R{{\mathbb R}}\def\PP{{\mathcal P}}
\def\Sym{{\rm Sym}} \def\D{{\mathbb D}} \def\d{{\bf d}}
\def\neg{{\rm neg}} \def\K{{\mathbb K}} \def\C{{\mathbb C}}
\def\HH{{\mathcal H}} \def\BB{{\mathcal B}} \def\UU{{\mathcal U}}
\def\CC{{\mathcal C}} \def\tr{{\rm tr}} \def\LL{{\mathcal L}}
\def\Z{{\mathbb Z}} \def\TR{{\rm TR}} \def\Syl{{\rm Syl}}
\def\DD{{\mathcal D}} \def\XX{{\mathcal X}} \def\Span{{\rm Span}}
\def\FF{{\mathcal F}} \def\EE{{\mathcal E}} \def\SS{{\mathcal S}}

\title{\bf{Singular Hecke algebras, Markov traces, and HOMFLY-type invariants}}
 
\author{\textsc{Luis Paris and Lo\"{\i}c Rabenda}}

\date{\today}

\maketitle

\begin{abstract}
\noindent
We define the singular Hecke algebra $\HH (SB_n)$ as the quotient of the singular braid monoid algebra 
$\C (q) [SB_n]$ by the Hecke relations $\sigma_k^2 = (q-1) \sigma_k +q$, $1 \le k\le n-1$. We observe 
that $\HH (SB_n)$ has a natural graduation, $\HH (SB_n)= \oplus_{d=0}^{+\infty} \HH (S_dB_n)$, where 
$\HH (S_dB_n)$ is the linear subspace spanned by the braids with $d$ singular points, and we show that 
$\HH (S_dB_n)$ is of finite dimension for all $d \ge 0$. We define a Markov trace on the sequence $\{ 
\HH( S_dB_n)\}_{n=1}^{+\infty}$ in the same way as for the Markov traces on the tower of 
(non-singular) Hecke algebras of the symmetric groups. We prove that a Markov trace determines an invariant 
on the links with a fixed number $d$ of singular points which satisfies some skein 
relation. Conversely, we prove that any invariant which satisfies this skein relation is of this form. 
Let $\TR_d$ denote the set of Markov traces on $\{ \HH (S_dB_n) \}_{n=1}^{+\infty}$. This is a $\C 
(q,z)$-vector space. Our main result is that $\TR_d$ is of dimension $d+1$. This result is completed 
with an explicit construction of a basis of $\TR_d$. Thanks to this result, we define a universal 
Markov trace and a universal HOMFLY-type invariant $\hat I: \LL \to \C[t^{\pm 1}, x^{\pm 1}, X,Y]$, where 
$\LL$ is the set of all (isotopy classes of) singular links, and $t,x,X,Y$ are variables. We show that 
this invariant is the unique invariant which satisfies some skein relation and some desingularization 
relation, and which takes the value 1 on the trivial knot. It is universal in the sense that, given two 
links $L,L'$ with $d$ singular points, we have $\hat I(L) = \hat I(L')$ if and only if $I(L) = I(L')$ 
for every invariant $I: \LL \to \C [^{\pm 1}, x^{\pm 1}]$ which satisfies the studied skein relation.
\end{abstract}

\noindent
{\bf AMS Subject Classification.}  Primary: 57M25. Secondary: 20C08, 20F36. 

\section{Introduction}

The {\it Hecke algebra} $\HH (B_n)$ of the symmetric group is a one parameter deformation of the 
symmetric group algebra studied in representation theory as well as in knot theory. Let $\K= \C(q)$ be 
the field of rational functions on a variable $q$, and let $B_n$ denote the braid group on $n$ 
strands. Then $\HH (B_n)$ is the quotient of the group algebra $\K [B_n]$ by the relations
\[
\sigma_k^2 = (q-1) \sigma_k +q\,, \quad 1 \le k\le n-1\,,
\]
where $\sigma_1, \dots, \sigma_{n-1}$ are the standard generators of $B_n$.

\bigskip\noindent
Let $z$ be a new variable. A {\it Markov trace} on the tower of algebras $\{ \HH 
(B_n)\}_{n=1}^{+\infty}$ is defined to be a collection of $\K$-linear maps
\[
\tr_n: \HH (B_n) \to \K(z)\,, \quad n\ge 1\,,
\]
such that
\begin{itemize}
\item
$\tr_n(\alpha \beta) = \tr_n (\beta \alpha)$ for all $\alpha, \beta \in B_n$ and all $n \ge 1$;
\item
$\tr_{n+1} (\beta) = \tr_n (\beta)$ for all $\beta \in B_n \subset B_{n+1}$ and all $n \ge 1$;
\item
$\tr_{n+1} (\beta \sigma_n) = z \cdot \tr_n (\beta)$ for all $\beta \in B_n$ and all $n \ge 1$.
\end{itemize}

\bigskip\noindent
Let $\LL_0$ denote the set of (isotopy classes of) links in $\R^3$. According to Jones \cite{Jones1}, a 
Markov trace $T= \{ \tr_n \}_{n=1}^{+\infty}$ on $\{ \HH (B_n) \}_{n=1}^{+\infty}$ determines a link 
invariant $I_T: \LL_0 \to \K (\sqrt{y})$, where $y=\frac{z-q+1}{qz}$. On the other hand, by a result of 
Ocneanu (see \cite{Jones1}, \cite{HOMFLY}), there exists a unique Markov trace which takes the value 1 
on the identity. In particular, the set of Markov traces form a one dimensional $\K(z)$-vector space 
spanned by the Ocneanu trace.

\bigskip\noindent
Let $A$ be an abelian group, let $I: \LL_0 \to A$ be an invariant, and let $t,x \in A$. We say that $I$ 
satisfies the {\it $(t,x)$ skein relation} if
\[
t^{-1} \cdot I(L_+) -t \cdot I(L_-) = x \cdot I(L_0)\,,
\]
for all links $L_+, L_-, L_0 \in \LL_0$ that have the same link diagram except in the neighborhood of a 
crossing where they are like in Figure 1.1. It is well-known that there exists a unique invariant $I: 
\LL_0 \to \C [t^{\pm 1}, x^{\pm 1}]$ which satisfies the $(t,x)$ skein relation and which takes the 
value 1 on the trivial knot. This invariant is equal to $I_T$ (up to a change of variables), 
where $T$ is the Ocneanu trace, and it is 
called the {\it HOMFLY polynomial} (see \cite{HOMFLY}, \cite{Jones2}, \cite{Jones1}, \cite{PrzTra1}).

\begin{figure}[htb]
\bigskip
\centerline{
\setlength{\unitlength}{.5cm}
\begin{picture}(15,4)
\put(0,1){\includegraphics[width=7.5cm]{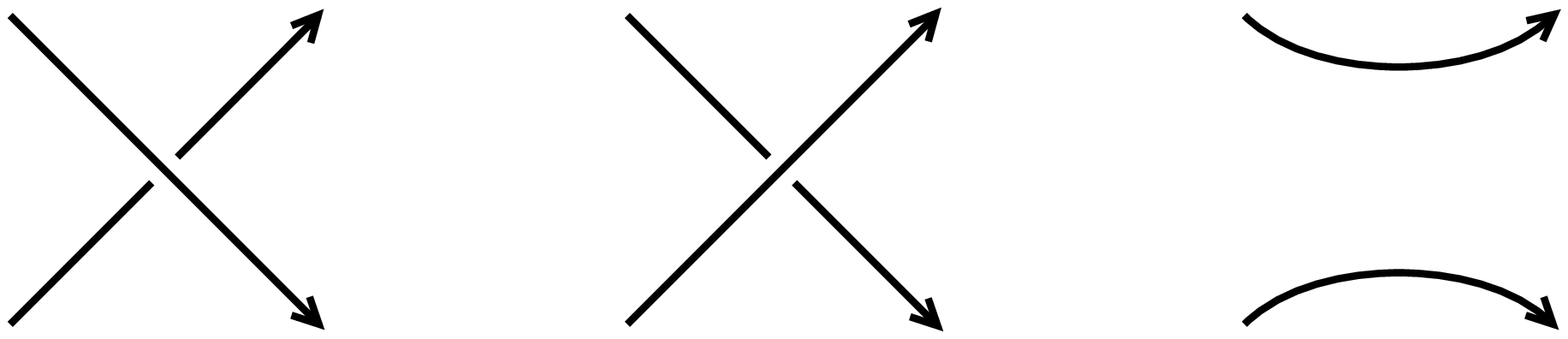}}
\put(1,0){\small $L_+$}
\put(7,0){\small $L_-$}
\put(13,0){\small $L_0$}
\end{picture}}
\bigskip
\centerline{{\bf Figure 1.1.} The links $L_+$, $L_-$, and $L_0$.}
\end{figure}

\bigskip\noindent
Our goal in this paper is to extend these constructions to the singular braids and links.

\bigskip\noindent
Let $SB_n$ denote the monoid of singular braids on $n$ strands. After some preliminaries on singular 
links and braids in Section 2, we develop the study of singular Hecke algebras, Markov traces, and 
related singular link invariants in Section 3. We define the {\it singular Hecke algebra} in a na\"{\i}ve 
way, as the quotient of the singular braid monoid algebra $\K [SB_n]$ by the Hecke relations
\[
\sigma_k^2 = (q-1) \sigma_k +q\,, \quad 1 \le k\le n-1\,.
\]

\bigskip\noindent
For $d\ge 0$, let $S_dB_n$ denote the set of braids with $d$ singular points. The algebra $\HH(SB_n)$ has a natural 
graduation
\[
\HH (SB_n) = \bigoplus_{d=0}^{+\infty} \HH (S_dB_n)\,,
\]
where $\HH (S_dB_n)$ is the linear subspace of $\HH (SB_n)$ spanned by $S_dB_n$. The algebra $\HH 
(SB_n)$ itself is of infinite dimension, but we show that each subspace $\HH (S_dB_n)$ of the 
graduation is of finite dimension over $\K$ (see Proposition 3.1).

\bigskip\noindent
A Markov trace on the sequence $\{ \HH (S_dB_n) \}_{n=1}^{+\infty}$ is defined in the same way as a 
Markov trace on $\{ \HH (B_n) \}_{n=1}^{+\infty}$. Let $\LL_d$ denote the set of (isotopy classes of) 
links with $d$ singular points. We prove that a Markov trace $T$ on $\{ \HH (S_dB_n) \}_{n=1}^{+\infty}$ 
determines an invariant $I_T: \LL_d \to \K (\sqrt{y})$ (see Proposition 3.3), and that this invariant 
satisfies the $(t,x)$ skein relation for $t= \sqrt{q} \sqrt{y}$ and $x= \sqrt{q} - \frac{1}{\sqrt{q}}$ 
(see Proposition 3.4). Conversely, any invariant $I: \LL_d \to \C (\sqrt{q}, \sqrt{y})$ which satisfies 
the $(t,x)$ skein relation is of the form $I=I_T$, where $T$ is a Markov trace on $\{ \HH (S_dB_n) 
\}_{n=1}^{+\infty}$ (with coefficients in $\C(\sqrt{q},\sqrt{y})$).

\bigskip\noindent
Section 4 contains the main result of the paper. Let $\TR_d$ denote the set of traces on $\HH 
(S_dB_n)$. This is a $\K(z)$-vector space. We prove that the dimension of $\TR_d$ is $d+1$, and construct an 
explicit basis of $\TR_d$ (see Theorem 4.7).

\bigskip\noindent
Let $\LL$ denote the set of all (isotopy classes of) singular links. Thanks to Section 4, we define in 
Section 5 a universal trace and a universal HOMFLY-type invariant $\hat I: \LL \to 
\C(\sqrt{q},\sqrt{y})[X,Y]$, where $y,X,Y$ are variables. We prove that $\hat I$ distinguishes two 
singular links $L,L' \in \LL_d$ (where $d$ is fixed) if and only if there exists an invariant $I: \LL_d 
\to \C (\sqrt{q}, \sqrt{y})$ which satisfies the $(t,x)$ skein relation for $t=\sqrt{q} 
\sqrt{y}$ and $x= \sqrt{q} - \frac{1}{\sqrt{q}}$, and which distinguishes $L$ and $L'$ (see Theorem~5.3). 
We also prove that $\hat I$ is the unique invariant with values in $\C [t^{\pm 1}, x^{\pm 1}, 
X,Y]$, which satisfies the $(t,x)$ skein relation and some desingularization relation, and which takes 
the value 1 on the trivial knot (see Proposition 5.4 and Theorem 5.5).

\bigskip\noindent
Our invariant $\hat I$ is more or less equivalent to the invariant of Kauffman and Vogel defined in 
\cite{KauVog1}. More precisely, the invariant of Kauffman and Vogel is the specialization $X=1$ of our 
invariant, but this specialization does not make much difference. Nevertheless, their approach is 
different from ours in the sense that they use singular Reidemeister moves to prove that their 
invariant is an invariant. They define some ``generalized Hecke algebras'' and define a Markov trace on 
this family of generalized Hecke algebras from which they can recover their invariant, but they 
definition of generalized Hecke algebras involves many relations besides the Hecke ones that are not 
natural in the context of an algebraic study. 

\bigskip\noindent
Note that the notion of Markov traces to study singular braids and links is also present in 
\cite{Baez1}, but the considered algebra in this paper is a one parameter deformation of $\C[B_n]$, 
which is specially adapted to the study Vassiliev invariants, but which has nothing to do with the Hecke 
relations. 
 

\section{Singular links and braids}

Let $n \ge 1$, and let $\S_1, \dots, \S_n$ be $n$ copies of the circle $\S^1$. A {\it singular link on 
$n$ components} is defined to be a smooth immersion $L: \S_1 \sqcup \cdots \sqcup \S_n \to \R^3$ such 
that, for all $P \in L$ (we identify the map $L$ with its image) there exist an open neighborhood $U$ 
of $P$ and a diffeomorphism $\varphi_U : U \to (-1,1)^3$ such that
\begin{itemize}
\item
$\varphi_U(P)= (0,0,0)$;
\item
either $\varphi_U (U \cap L) = (-1,1) \times \{0\} \times \{0\}$, or $\varphi_U (U \cap L)= ((-1,1) 
\times \{0\} \times \{0\}) \cup (\{0\} \times (-1,1) \times \{0\})$.
\end{itemize}
In the second case, when $\varphi_U (U \cap L)= ((-1,1) \times \{0\} \times \{0\}) \cup (\{0\} \times 
(-1,1) \times \{0\})$, we say that $P$ is a {\it singular point} of $L$. 

\bigskip\noindent
In this context, the 
admissible isotopies preserve the whole structure. Moreover, the circle $\S^1$ as well as the links are 
always assumed to be oriented.

\bigskip\noindent
Let $\pi: \R^3 \to \R^2$, $(x,y,z) \mapsto (x,y)$, be the projection on the two first coordinates. Up to 
isotopy, we can assume that $\pi \circ L: \S_1 \sqcup \cdots \sqcup \S_n \to \R^2$ is a smooth 
immersion. Moreover, for all $P \in \pi (L)$, we can assume that there exist a neighborhood $V$ of $P$ 
and a diffeomorphism $\psi_V: V \to (-1,1)^2$ such that
\begin{itemize}
\item
$\psi_V (P) = (0,0)$;
\item
either $\psi_V (V \cap \pi(L)) = (-1,1) \times \{0\}$, or $\psi_V (V \cap \pi(L)) = ((-1,1) \times 
\{0\}) \cup (\{0\} \times (-1,1))$.
\end{itemize}
In the second case, when $\psi_V (V \cap \pi(L)) = ((-1,1) \times \{0\}) \cup (\{0\} \times (-1,1))$, we 
say that $P$ is a {\it crossing} of $\pi(L)$, and we indicate graphically like in Figure 2.1 if $P$ is 
the preimage of a singular point of $L$, or else which strand passes over the other. Such a 
graphical representation of $L$ is called a {\it link diagram} of $L$.

\begin{figure}[htb]
\bigskip
\centerline{
\setlength{\unitlength}{.5cm}
\begin{picture}(15,5)
\put(0,2){\includegraphics[width=7.5cm]{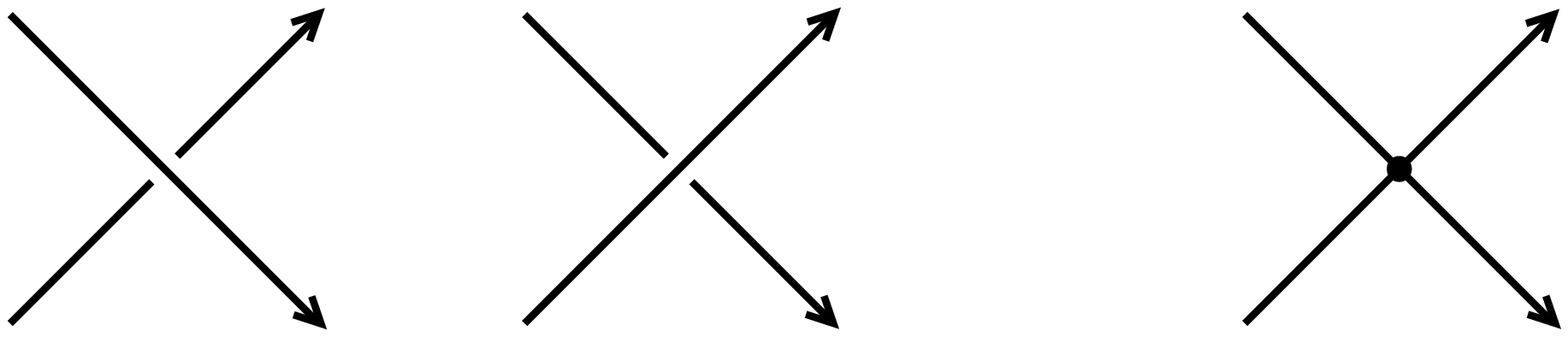}}
\put(2.5,1){\small non-singular}
\put(3,0){\small crossings}
\put(12.3,1){\small singular}
\put(12.3,0){\small crossing}
\end{picture}}
\bigskip
\centerline{{\bf Figure 2.1.}  Crossings in a singular link diagram.}
\end{figure}

\bigskip\noindent
Obviously, a singular link can be determined by many link diagrams. However, we have the following.

\bigskip\noindent
{\bf Theorem 2.1} (Kauffman \cite{Kauff1}). {\it Two link diagrams represent the same singular link up 
to isotopy if and only if one can pass from one to the other by a finite sequence of singular 
Reidemeister moves as shown in Figure 2.2.}
\qed

\begin{figure}[htb]
\bigskip
\centerline{
\setlength{\unitlength}{.5cm}
\begin{picture}(30,22)
\put(0,0){\includegraphics[width=15cm]{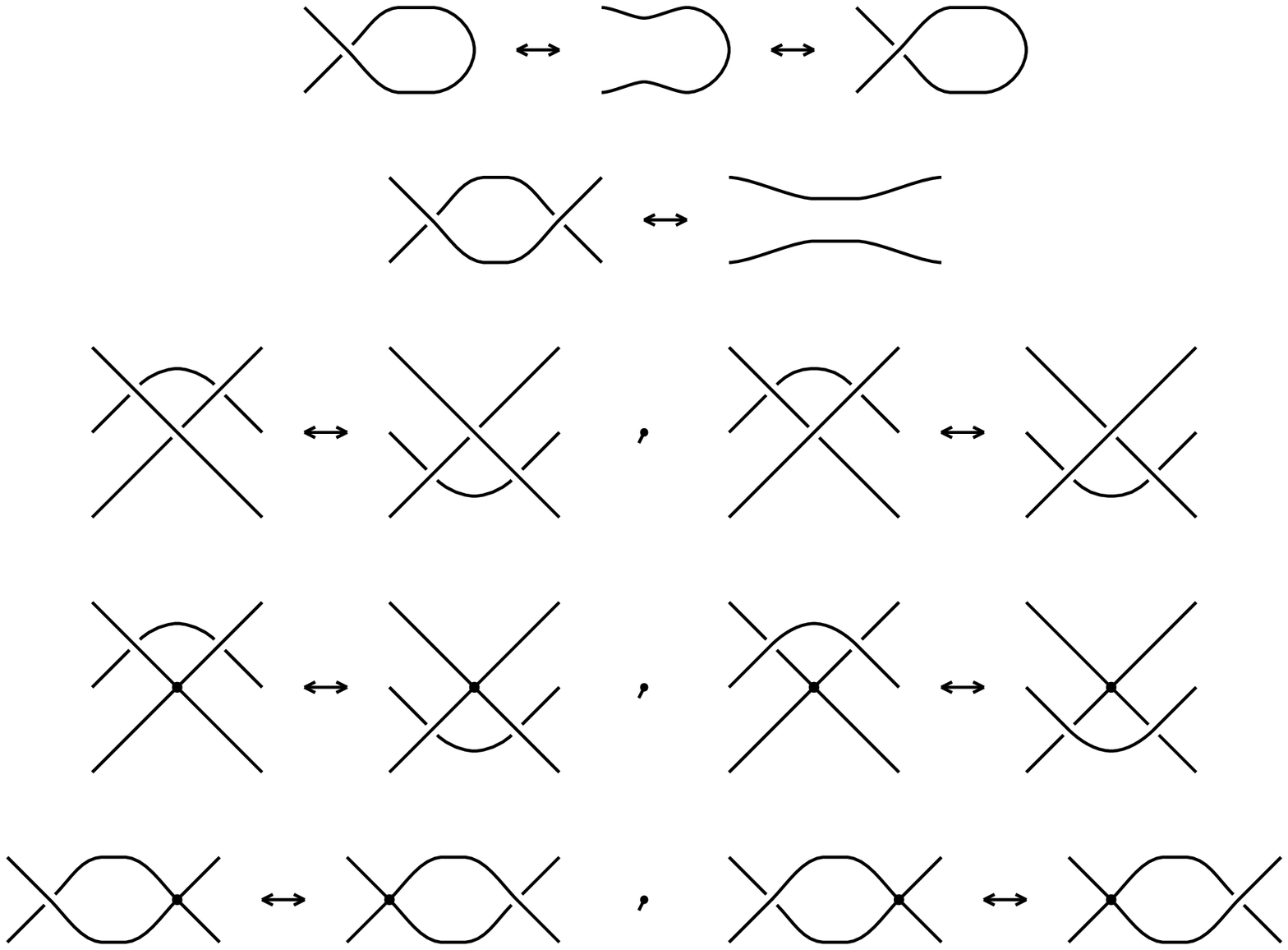}}
\end{picture}}
\bigskip
\centerline{{\bf Figure 2.2.} Singular Reidemeister moves.}
\end{figure}

\bigskip\noindent
Let $\PP = \{ P_1, \dots, P_n\}$ be a set of $n$ punctures in $\R^2$ (except mention of the 
contrary, we will always assume $P_k = (k,0)$ for all $1 \le k\le n$). A {\it singular braid on $n$ strands 
based at $\PP$} is defined to be a $n$-tuple $\beta= (b_1, \dots, b_n)$ of smooth paths, $b_k: [0,1] 
\to \R^2 \times [0,1]$, such that:
\begin{itemize}
\item
There exists a permutation $\chi \in \Sym_n$ such that $b_k (0)= (P_k, 0)$ and $b_k (1) = (P_{\chi(k)}, 
1)$ for all $1 \le k\le n$.
\item
Let $\pi_3: \R^2 \times [0,1] \to [0,1]$, $(x,y,t) \mapsto t$, be the projection on the third 
coordinate. Then $\pi_3 (b_k(t)) = t$ for all $1 \le k \le n$ and all $t \in [0,1]$.
\item
Let
\[
\breve \beta = b_1 ((0,1)) \cup \cdots \cup b_n ((0,1)) \subset \R^2 \times (0,1)\,.
\]
For all $P \in \breve \beta$, there exist a neighborhood $U$ of $P$ and a diffeomorphism $\varphi_U: U 
\to (-1,1)^3$ such that $\varphi_U (P) = (0,0,0)$, and either $\varphi_U( U \cap \breve \beta) = (-1,1) 
\times \{0\} \times \{0\}$, or $\varphi_U (U \cap \breve \beta)= ((-1,1) \times \{0\} \times \{0\}) 
\cup (\{0\} \times (-1,1) \times \{0\})$.
\end{itemize}
In the above definition, we call the point $P \in \breve\beta$ a {\it singular point} of $\beta$ if $\varphi_U (U \cap 
\breve \beta)= ((-1,1) \times \{0\} \times \{0\}) \cup (\{0\} \times (-1,1) \times \{0\})$.

\bigskip\noindent
The isotopy classes of singular braids form a monoid (and not a group) called the {\it singular braid 
monoid on $n$ strands} and denoted by $SB_n$. The monoid operation is the concatenation.

\bigskip\noindent
Let $\pi: \R^2 \times [0,1] \to \R \times [0,1]$, $(x,y,t) \mapsto (x,t)$. Up to isotopy, we can assume 
that, for all $P \in \pi (\breve \beta)$, there exist a neighborhood $V$ of $P$ and a diffeomorphism 
$\psi_V: V \to (-1,1)^2$ such that $\psi_V(P)= (0,0)$, and either $\psi_V (V \cap \pi (\breve \beta)) = 
(-1,1) \times \{0\}$, or $\psi_V (V \cap \pi (\breve \beta)) = ((-1,1) \times \{0\}) \cup (\{0\} \times 
(-1,1))$. In the second case, when $\psi_V (V \cap \pi (\breve \beta)) = ((-1,1) \times \{0\}) \cup 
(\{0\} \times (-1,1))$, we say that $P$ is a {\it crossing} and we indicate graphically like in Figure 
2.1 if $P$ is the preimage of a singular point of $\beta$, or else which strand passes over the other. 
Such a representation of $\beta$ is called a {\it braid diagram} of $\beta$.

\bigskip\noindent
{\bf Theorem 2.2} (Baez \cite{Baez1}, Birman \cite{Birma1}). {\it The monoid $SB_n$ has a monoid 
presentation with generators
\[
\sigma_1, \dots, , \sigma_{n-1}, \sigma_1^{-1}, \dots, \sigma_{n-1}^{-1}, \tau_1, \dots, \tau_{n-1}\,,
\]
and relations
\[\begin{array}{cl}
\sigma_k \sigma_k^{-1} = \sigma_k^{-1} \sigma_k = 1 &\quad\text{for } 1 \le k\le n-1\,,\\
\sigma_k \tau_k = \tau_k \sigma_k &\quad \text{for } 1 \le k\le n-1\,,\\
\sigma_k \sigma_l \sigma_k = \sigma_l \sigma_k \sigma_l &\quad\text{if } |k-l| =1\,,\\
\sigma_k \sigma_l \tau_k = \tau_l \sigma_k \sigma_l &\quad\text{if }|k-l|=1\,,\\
\sigma_k \sigma_l = \sigma_l \sigma_k &\quad\text{if } |k-l|\ge 2\,,\\
\sigma_k \tau_l = \tau_l \sigma_k &\quad{if }|k-l|\ge 2\,,\\
\tau_k \tau_l = \tau_l \tau_k &\quad\text{if } |k-l|\ge 2\,.
\end{array}\]}
\qed

\bigskip\noindent
The braid $\sigma_k$ in the above presentation is the standard $k$-th generator of the braid group $B_n$ 
(see Figure 2.3). The braid $\tau_k$ is a singular braid with a single singular point which involves 
the $k$-th strand and the $(k+1)$-th strand (see Figure 2.3).

\begin{figure}[htb]
\bigskip
\centerline{
\setlength{\unitlength}{.5cm}
\begin{picture}(18,7)
\put(4,0){\includegraphics[width=7cm]{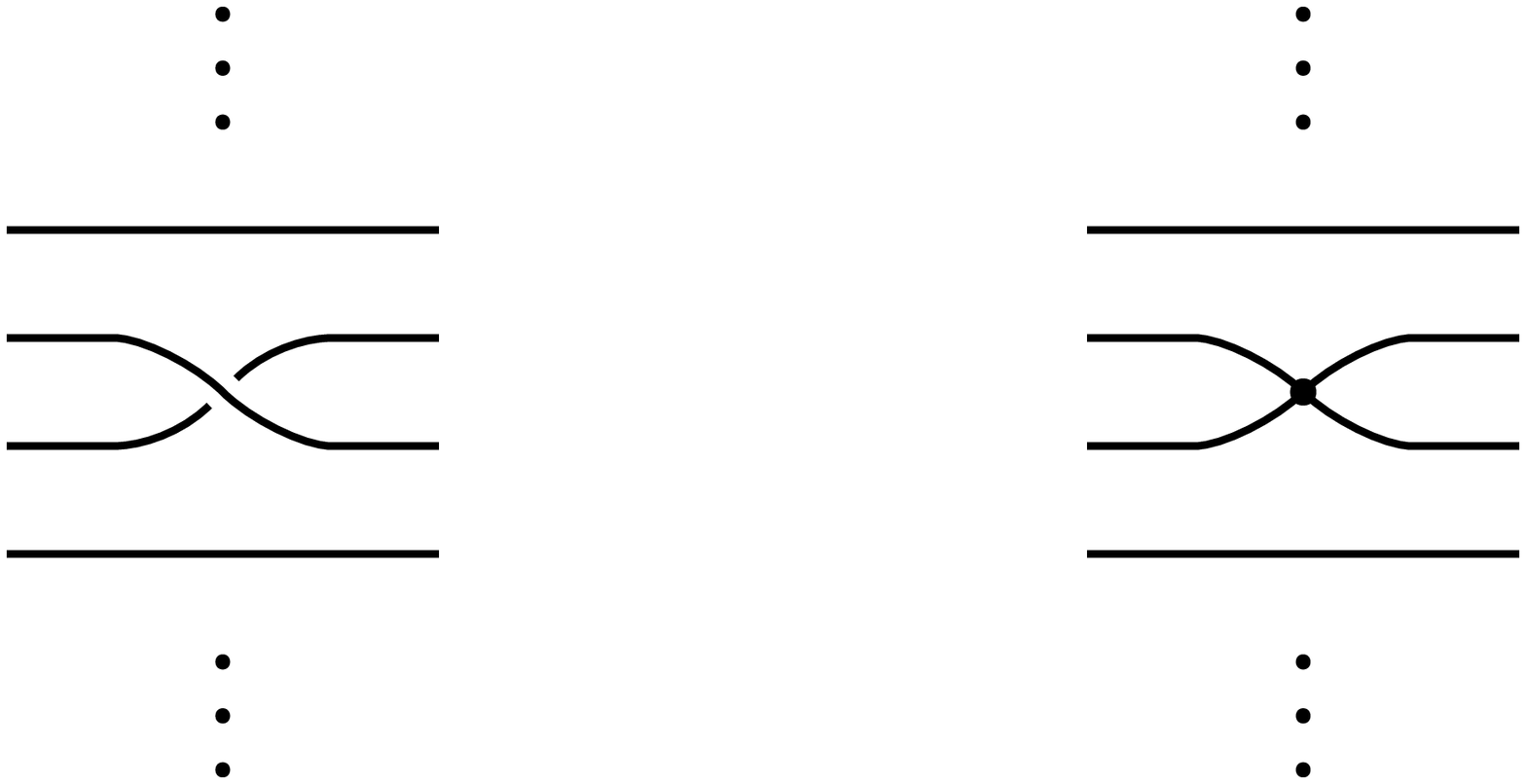}}
\put(0.5,3){$\sigma_k=$}
\put(3.3,3){\small $k$}
\put(2.3,4){\small $k+1$}
\put(10.5,3){$\tau_k=$}
\put(13.3,3){\small $k$}
\put(12.3,4){\small $k+1$}
\end{picture}}
\bigskip
\centerline{{\bf Figure 2.3.} Generators of $SB_n$.}
\end{figure}

\bigskip\noindent
From a singular braid $\beta= (b_1, \dots, b_n)$ one can construct a singular link, called the {\it 
closure} of $\beta$ and denoted by $\hat \beta$, as follows. Let $\varphi_1 : \R^2 \times [0,1] \to 
\R^2 \times \S^1$ be the map defined by $\varphi_1( P,t)= (P, e^{2i\pi t})$ for all $(P,t) \in \R^2 
\times [0,1]$. In other words, $\varphi_1$ identifies $(P,0)$ with $(P,1)$ for all $P \in \R^2$. Let 
$\breve \D$ denote the interior of the unit disk. We take a diffeomorphism $\psi: \R^2 \to \breve \D$ 
and we extend it to a diffeomorphism $\varphi_2: \R^2 \times \S^1 \to \breve \D \times \S^1$, $(P, \xi) 
\mapsto (\psi(P), \xi)$. For instance, the map $\psi$ can be defined by
\[
\psi(x,y)= \left( {x \over 1+ \sqrt{x^2+y^2}}, {y  \over 1+ \sqrt{x^2+y^2}} \right) \quad \text{for all 
} (x,y) \in \R^2\,.
\]
Finally, we take a standard embedding $\varphi_3: \breve \D \times \S^1 \to \R^3$ and we set $\hat \beta 
= (\varphi_3 \circ \varphi_2 \circ \varphi_1) (\beta)$. For instance, the map $\varphi_3$ can be defined 
by
\[
\varphi_3 (x,y,e^{2i\pi t}) = ((2+x) \cos 2\pi t, (2+x) \sin 2\pi t, y)\,.
\]
An example of a closed singular braid is illustrated in Figure 2.4.

\begin{figure}[htb]
\bigskip
\centerline{
\setlength{\unitlength}{.5cm}
\begin{picture}(20.5,9)
\put(2,0){\includegraphics[width=9.25cm]{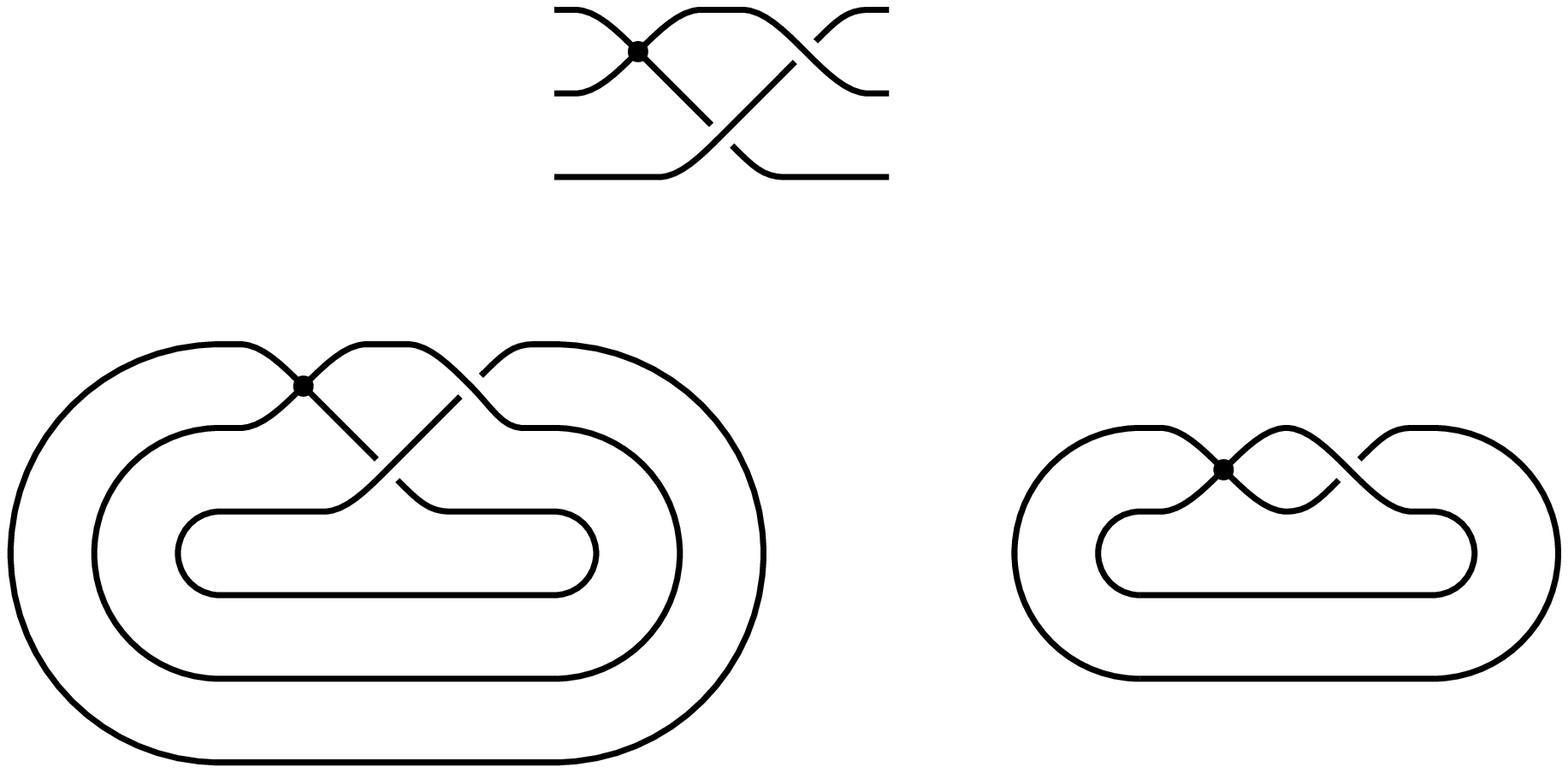}}
\put(0,2){$\hat\beta =$}
\put(12,2){$\simeq$}
\put(6.5,8){$\beta=$}
\end{picture}}
\bigskip
\centerline{{\bf Figure 2.4.} A closed singular braid.}
\end{figure}

\bigskip\noindent
The following theorem is known as the Alexander theorem for singular links, and can be proved with 
almost the same arguments as for the classical Alexander theorem. We give a proof here because it will be 
needed in the next sections to get the skein relation for our invariants. This proof is based on the 
sketch given by Birman in \cite{Birma1}.

\bigskip\noindent
{\bf Theorem 2.3} (Birman \cite{Birma1}). {\it Every singular link is a closed braid.}

\bigskip\noindent
{\bf Proof.} Let $L$ be a singular link. We consider the projection $\pi: \R^3 \to \R^2$, $(x,y,z) 
\mapsto (x,y)$ and we assume that $\pi (L)$ is a link diagram. We also assume that $(0,0) \not\in \pi 
(L)$, so that the vertical line $\d = \{ (0,0,z) ; z \in \R \}$ does not touch $L$.

\bigskip\noindent
Up to isotopy, we assume that there exist intervals $I_1, \dots, I_p$ in $\R^3$ such that
\begin{enumerate}
\item
$L= I_1 \cup I_2 \cup \cdots \cup I_p$;
\item
$I_i \cap I_j$ is either empty or a single point which must be an extremity of both, $I_i$ and $I_j$, 
for all $1 \le i \neq j \le p$;
\item
$I_i$ and $\d$ are not coplanar for any $1 \le i\le p$;
\item
no extremity of $\pi(I_i)$ is a crossing of $\pi (L)$ for all $1 \le i \le p$.
\end{enumerate}
Note that the above assumptions imply that $L$ is piecewise linear but not smooth anymore. However, 
we can smooth the corners along the whole proof.

\bigskip\noindent
Let $P$ be a crossing of $\pi (L)$, and let $I_i$ and $I_j$ be the intervals such that $P \in \pi(I_i) 
\cap \pi(I_j)$. Let $P_i= (P,z_i) \in I_i$ (resp. $P_j= (P,z_j)\in I_j$) such that $\pi(P_i)=P$ (resp. 
$\pi(P_j)=P$). We say that $P$ is an {\it upper crossing} of $I_i$ if $z_i > z_j$, that $P$ is a {\it 
lower crossing} of $I_i$ if $z_i<z_j$, and that $P$ is a {\it singular crossing} of $I_i$ if $z_i=z_j$. 
We say that an interval $I_i$ is an {\it upper interval} (resp. {\it lower interval}, {\it singular 
interval}) if all the crossings on $I_i$ are upper crossings (resp. lower crossings, singular 
crossings) of $I_i$. We use the convention that $I_i$ is an upper interval if there is no crossing on 
it.

\bigskip\noindent
Without loss of generality, we can add the following assumption on the intervals.
\begin{enumerate}
\setcounter{enumi}{4}
\item
Every interval is either upper, or lower, or singular.
\end{enumerate}
Recall that $L$ is oriented, and this orientation induces an orientation on each $I_i$. We set 
$\pi(I_i)= [A_i,B_i]$ so that the orientation goes from $A_i$ to $B_i$. Not that the hypothesis (3) 
($I_i$ and $\d$ are not coplanar) implies that $A_i \neq B_i$ and $O=(0,0)$ does not lie in the line 
$\langle A_i, B_i \rangle$ spanned by $A_i$ and $B_i$. We endow $\R^2$ with the standard orientation, 
and we say that $I_i$ is {\it positive} if $( \overrightarrow{ A_i B_i}, \overrightarrow{A_iO})$ is a 
direct basis of $\R^2$ (see Figure 2.5), and we say that $I_i$ is {\it negative} otherwise.

\begin{figure}[htb]
\bigskip
\centerline{
\setlength{\unitlength}{.5cm}
\begin{picture}(20,8)
\put(0,0){\includegraphics[width=10cm]{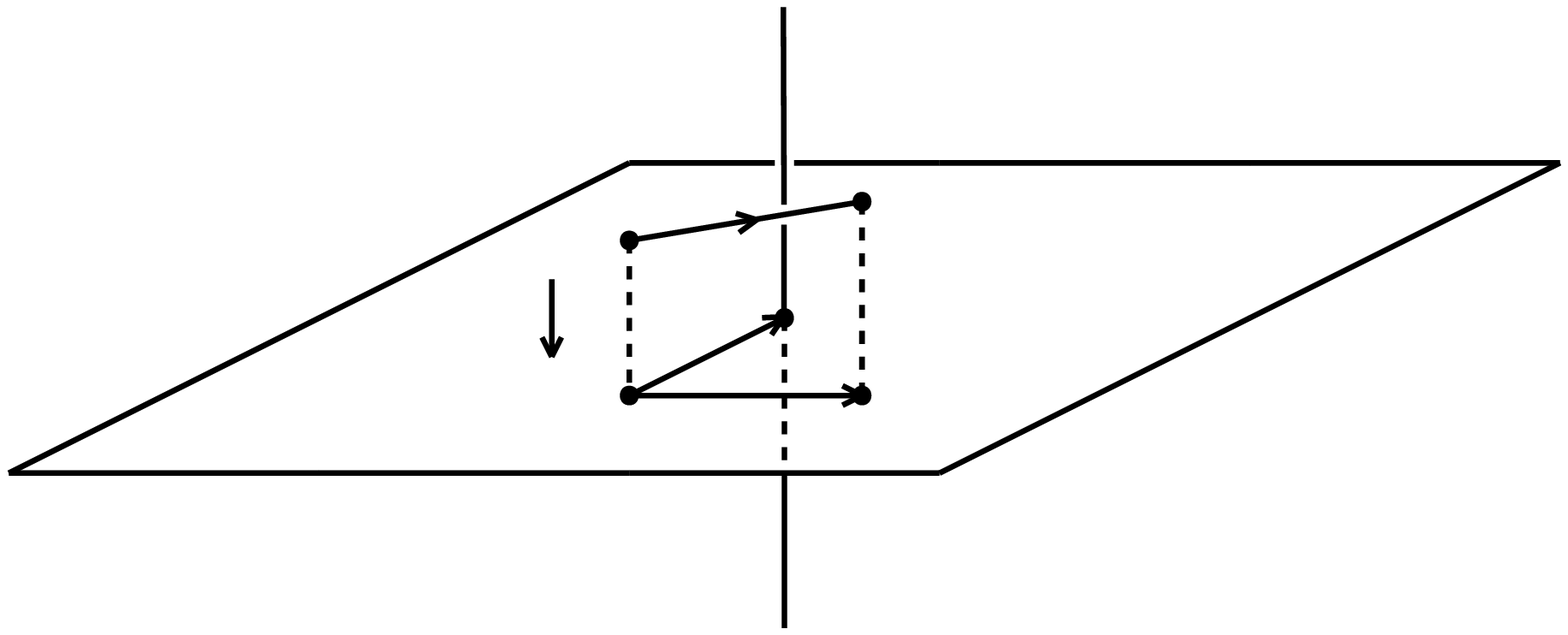}}
\put(7.3,2.5){\small $A_i$}
\put(11.3,2.5){\small $B_i$}
\put(10.2,4){\small $O$}
\put(6.3,3.9){\small $\pi$}
\put(8.5,5.4){\small $I_i$}
\put(10.4,6.5){\small ${\bf d}$}
\end{picture}}
\bigskip
\centerline{{\bf Figure 2.5.} A positive interval.}
\end{figure}

\bigskip\noindent
Let $\neg (L)$ denote the number of negative intervals. It is easily seen that $L$ is a closed braid 
if $\neg (L) = 0$.

\bigskip\noindent
Upon transforming $L$ and $\pi (L)$ and adding new (non-singular) crossings, we can also assume that
\begin{enumerate}
\setcounter{enumi}{5}
\item
Every singular interval is positive.
\end{enumerate}
From now on in the proof, all the isotopies we shall consider preserve Conditions (1) to (6). We prove 
by induction on $\neg (L)$ that $L$ is a closed braid.

\bigskip\noindent
Suppose $\neg(L) >0$. Let $I_i$ be a negative interval. By Condition (6), $I_i$ is either upper or 
lower (say it is upper). It is easily seen that there exists a flat triangle $T$ embedded in $\R^3$ 
such that $T \cap L = I_i$, and $T \cap \d$ is a unique point in the interior of $T$ (see Figure 2.6). 
Let $I_{i\,1}$ and $I_{i\,2}$ be the other sides of $T$, and set
\[
L'= (L \setminus I_i) \cup (I_{i\,1} \cup I_{i\,2})\,.
\]
Clearly, we can choose $T$ so that $I_{i\,1}$ and $I_{i\,2}$ are both upper and positive intervals. Then 
$L'$ is isotopic to $L$, it satisfies Conditions (1) to (6), and $\neg (L') = \neg(L)-1$.
\qed

\begin{figure}[htb]
\bigskip
\centerline{
\setlength{\unitlength}{.5cm}
\begin{picture}(5,5)
\put(1,0){\includegraphics[width=2cm]{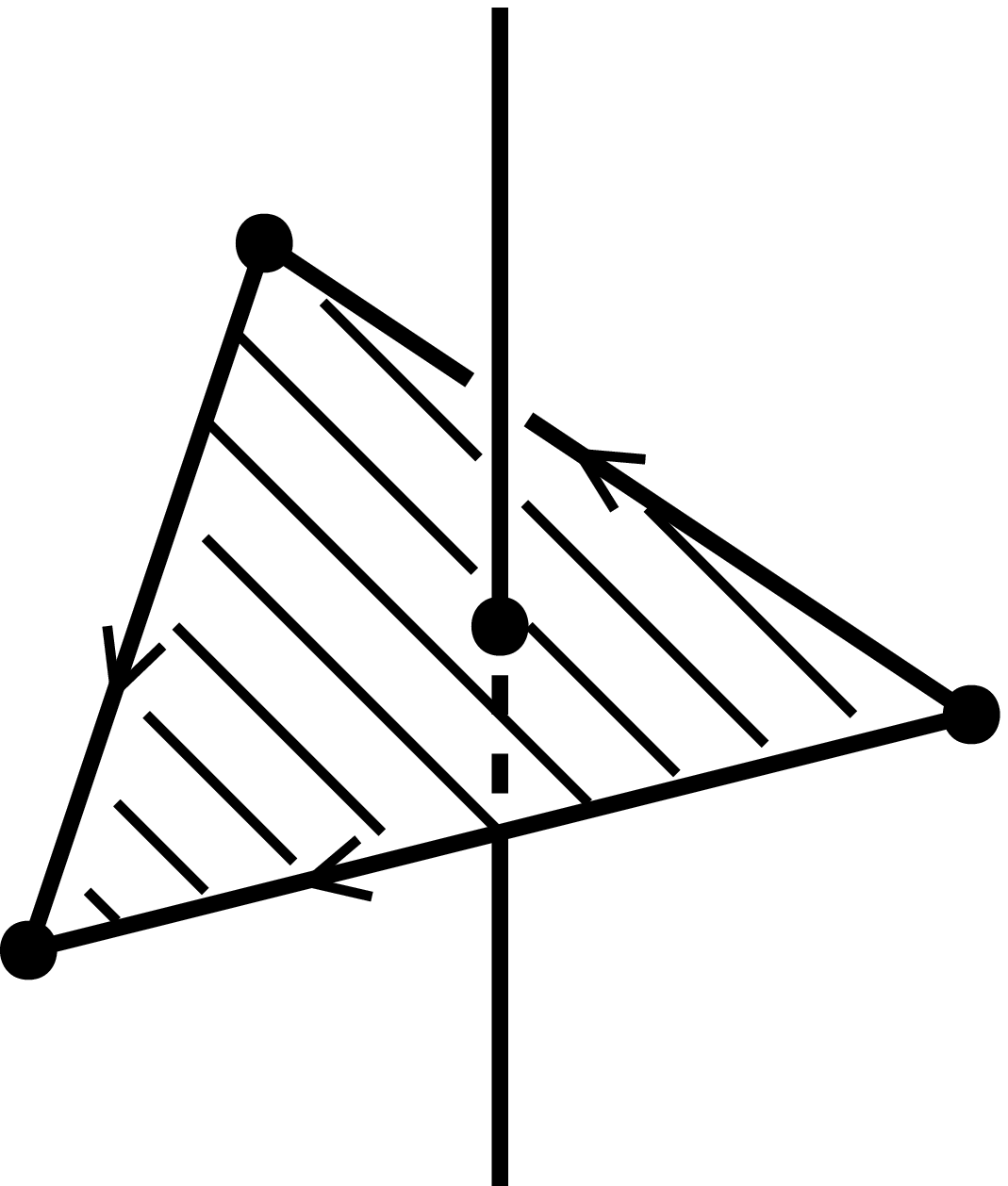}}
\put(1.8,0.5){\small $I_i$}
\put(0.4,2.5){\small $I_{i\,2}$}
\put(3.8,2.8){\small $I_{i\,1}$}
\end{picture}}
\bigskip
\centerline{{\bf Figure 2.6.} A triangle.}
\end{figure}

\bigskip\noindent
Now, consider the set $\sqcup_{n=1}^{+\infty} SB_n$ of all singular braids. We may often use the notation 
$(\beta,n)$ to denote a braid $\beta \in SB_n$ in case we want to emphasize the number $n$ of strands.

\bigskip\noindent
We say that two singular braids $(\alpha, n)$ and $(\beta, m)$ are connected by a {\it Markov move} if 
either
\begin{itemize}
\item
$n=m$, $\alpha= \gamma_1 \gamma_2$, and $\beta = \gamma_2 \gamma_1$, for some $\gamma_1, \gamma_2 \in 
SB_n$; or
\item
$n=m+1$ and $\alpha= \beta \sigma_n^{\pm 1}$; or
\item
$m=n+1$ and $\beta= \alpha \sigma_m^{\pm 1}$.
\end{itemize}

\bigskip\noindent
{\bf Theorem 2.4} (Gemein \cite{Gemei1}). {\it Let $(\alpha,n)$ and $(\beta,m)$ be two singular braids. Then 
$\hat \alpha$ and $\hat \beta$ are isotopic if and only if $(\alpha,n)$ and $(\beta,m)$ are connected 
by a finite sequence of Markov moves.}
\qed
 

\section{Singular Hecke algebras, Markov traces, and singular link invariants}

Recall that $\K = \C (q)$ denotes the field of rational functions on a variable $q$. We define the {\it singular 
Hecke algebra} $\HH (SB_n)$ to be the quotient of the monoid algebra $\K [SB_n]$ by the relations
\begin{equation}\label{R31}
\sigma_k^2 = (q-1) \sigma_k +q\,, \quad 1 \le k \le n-1\,.
\end{equation}

\bigskip\noindent
For $d \ge 0$, we denote by $S_dB_n$ the set of (isotopy classes of) singular braids on $n$ strands 
with $d$ singular points, and by $\K [S_dB_n]$ the $\K$-linear subspace of $\K [SB_n]$ spanned by 
$S_dB_n$. Note that $S_0B_n=B_n$ is the braid group, and $\K[S_0B_n] = \K [B_n]$ is the group algebra of 
$B_n$. We have the graduation
\[
\K [SB_n] = \bigoplus_{d=0}^{+\infty} \K [S_dB_n]\,.
\]
The relations \eqref{R31} that define the singular Hecke algebra involve only elements of degree zero, thus the 
graduation of $\K [SB_n]$ induces a graduation of $\HH (SB_n)$:
\[
\HH (SB_n) = \bigoplus_{d=0}^{+\infty} \HH (S_dB_n)\,,
\]
where $\HH (S_dB_n)$ is the $\K$-linear subspace of $\HH (SB_n)$ spanned by $S_dB_n$.

\bigskip\noindent
It is known that $\HH (B_n)$ has dimension $n!$, and has a basis $\BB_n$ which can be described as 
follows (see \cite{Jones1}). For $n \ge 2$ we set
\[
\UU_n = \{ 1, \sigma_{n-1}, \sigma_{n-1} \sigma_{n-2}, \dots, \sigma_{n-1} \sigma_{n-2} \cdots \sigma_2 
\sigma_1 \} \,.
\]
Then $\BB_n$ is defined by induction on $n$ by
\[
\BB_1 = \{1\}\,, \quad \BB_n = \{ \beta u\ ;\ \beta \in \BB_{n-1} \text{ and } u \in \UU_n\} \quad \text{ 
if } n \ge 2\,.
\]
The singular Hecke algebra $\HH (SB_n)$ is not of finite dimension, but each subspace $\HH 
(S_dB_n)$ of the graduation is of finite dimension. Indeed:

\bigskip\noindent
{\bf Proposition 3.1.} {\it Let $d \ge 0$, and let $n \ge 2$. Let $\CC_{d,n}$ denote the set of singular 
braids of the form $\tau_{i_1} \cdots \tau_{i_d} \beta$, where $1 \le i_j \le n-1$ for $1 \le j \le d$, 
and $\beta \in \BB_n$. Then $\CC_{d,n}$ spans $\HH (S_dB_n)$.}

\bigskip\noindent
{\bf Proof.} Observe that the Hecke relation \eqref{R31} implies that
\[
\sigma_k^{-1} = q^{-1} \sigma_k -q^{-1} (q-1)\,, \quad \text{for all } 1 \le k \le n-1\,.
\]
Let $i,j \in \{ 1, \dots, n-1 \}$ such that $|i-j|=1$, and let $a \ge 1$. We calculate $\sigma_i^2 
\tau_j^a$ in two ways. Firstly,
\[\begin{array}{rcl}
\sigma_i^2 \tau_j^a &=& \sigma_i \sigma_j^{-1} \sigma_j \sigma_i \tau_j^a\\
&=&\sigma_i \sigma_j^{-1} \tau_i^a \sigma_j \sigma_i\\
&=& q^{-1} \sigma_i \sigma_j \tau_i^a \sigma_j \sigma_i - q^{-1} (q-1) \sigma_i \tau_i^a \sigma_j 
\sigma_i\\
&=& q^{-1} \tau_j^a \sigma_i \sigma_j^2 \sigma_i - q^{-1} (q-1) \tau_i^a \sigma_i \sigma_j \sigma_i\\
&=& q^{-1} (q-1) \tau_j^a \sigma_i \sigma_j \sigma_i + (q-1) \tau_j^a \sigma_i + q \tau_j^a - q^{-1} 
(q-1) \tau_i^a \sigma_i \sigma_j \sigma_i
\end{array}\] 
Secondly,
\[
\sigma_i^2 \tau_j^a = (q-1) \sigma_i \tau_j^a + q \tau_j^a\,.
\]
These two equalities imply
\begin{equation}\label{R32}
\sigma_i \tau_j^a = q^{-1} \tau_j^a \sigma_i \sigma_j \sigma_i - q^{-1} \tau_i^a \sigma_i \sigma_j 
\sigma_i + \tau_j^a \sigma_i \,.
\end{equation}
On the other hand, by Theorem 2.2, if $i,j \in \{ 1, \dots, n-1\}$ are such that $|i-j| \neq 1$, then
\begin{equation}\label{R33}
\sigma_i \tau_j^a = \tau_j^a \sigma_i\,.
\end{equation}
The equalities \eqref{R32} and \eqref{R33} show that every element of $\HH (S_d B_n)$ is a linear 
combination of elements of the form $\tau_{i_1} \cdots \tau_{i_d} \omega$, where $1 \le i_j \le n-1$ 
for $1 \le j \le d$, and $\omega \in \HH (B_n)$. Now, since $\BB_n$ is a basis of $\HH(B_n)$, we 
conclude that every element of $\HH (S_dB_n)$ is a linear combination of elements of the form 
$\tau_{i_1} \cdots \tau_{i_d} \beta$, where $1 \le i_j \le n-1$ for $1 \le j \le d$, and $\beta \in 
\BB_n$.
\qed

\bigskip\noindent
However, $\CC_{d,n}$ is not a basis of $\HH (S_dB_n)$ in general. Indeed:

\bigskip\noindent
{\bf Lemma 3.2.} {\it Let $i,j \in \{ 1, \dots, n-1\}$ such that $|i-j|=1$, and let $a \ge 1$. Then
\begin{equation}\label{R34}
\begin{split}
&\tau_i^a (\sigma_i \sigma_j + \sigma_j \sigma_i - (q-1) \sigma_i - (q-1) \sigma_j + (q^2-q+1))\\ 
=\ &\tau_j^a (\sigma_i \sigma_j + \sigma_j \sigma_i - (q-1) \sigma_i - (q-1) \sigma_j + (q^2-q+1))\,.
\end{split}
\end{equation}} 

\bigskip\noindent
{\bf Proof.} Recall the equality \eqref{R32} in the proof of Proposition 3.1:
\[
\sigma_i \tau_j^a = q^{-1} \tau_j^a \sigma_i \sigma_j \sigma_i - q^{-1} \tau_i^a \sigma_i \sigma_j 
\sigma_i + \tau_j^a \sigma_i \,.
\]
We multiply this equality on the right hand side by $\sigma_i^{-1} \sigma_j^{-1}$ and we get
\[\begin{array}{crcl}
&\sigma_i \tau_j^a \sigma_i^{-1} \sigma_j^{-1} &=& q^{-1} \tau_j^a \sigma_i - q^{-1} \tau_i^a \sigma_i 
+ \tau_j^a \sigma_j^{-1}\\
\Leftrightarrow\quad &\sigma_i \sigma_i^{-1} \sigma_j^{-1}\tau_i^a  &=& q^{-1} \tau_j^a \sigma_i - q^{-
1} \tau_i^a \sigma_i + q^{-1}\tau_j^a \sigma_j - q^{-1} (q-1) \tau_j^a\\
\Leftrightarrow\quad & q^{-1}\sigma_j\tau_i^a - q^{-1} (q-1) \tau_i^a &=& q^{-1} \tau_j^a \sigma_i - 
q^{-1} \tau_i^a \sigma_i + q^{-1}\tau_j^a \sigma_j - q^{-1} (q-1) \tau_j^a
\end{array}\]
thus
\begin{equation}\label{R35} 
\sigma_j \tau_i^a = \tau_j^a (\sigma_i + \sigma_j - (q-1)) - \tau_i^a (\sigma_i - (q-1))\,.
\end{equation}
Now, we apply twice \eqref{R35} to $\sigma_i \sigma_j \tau_i^a$ and obtain
\[\begin{array}{rcl}
\sigma_i \sigma_j \tau_i^a &=& \sigma_i \tau_j^a (\sigma_i + \sigma_j -(q-1)) - \sigma_i \tau_i^a 
(\sigma_i -(q-1))\\
&=& \tau_i^a (\sigma_i + \sigma_j - (q-1))^2 - \tau_j^a (\sigma_j - (q-1)) (\sigma_i + \sigma_j - (q-
1)) - \tau_i^a \sigma_i (\sigma_i - (q-1))\\
&=& \tau_i^a (\sigma_i \sigma_j + \sigma_j \sigma_i - (q-1) \sigma_i - (q-1) \sigma_j + (q^2-q+1))\\ 
&&\hfill  - \tau_j^a (\sigma_j \sigma_i - (q-1) \sigma_i - (q-1) \sigma_j + (q^2-q+1))\,.\\
\end{array}\]
Since $\sigma_i \sigma_j \tau_i^a = \tau_j^a \sigma_i \sigma_j$, it follows that
\[
\begin{split}
&\tau_i^a (\sigma_i \sigma_j + \sigma_j \sigma_i - (q-1) \sigma_i - (q-1) \sigma_j + (q^2-q+1))\\ 
=\ &\tau_j^a (\sigma_i \sigma_j + \sigma_j \sigma_i - (q-1) \sigma_i - (q-1) \sigma_j + (q^2-q+1))\,.
\end{split}
\]
\qed

\bigskip\noindent
{\bf Remark.} We do not know the dimension of $\HH(S_dB_n)$ if $d \ge 1$ and $n \ge 3$.

\bigskip\noindent
We turn now to the definition of a Markov trace, but, before, we make the following remark.

\bigskip\noindent
{\bf Remark.}
The basis $\BB_n$ of $\HH (B_n)$ can be viewed as a subset of $\BB_{n+1}$. This implies that the 
natural embedding $B_n \hookrightarrow B_{n+1}$ leads to an injective homomorphism $\HH (B_n) 
\hookrightarrow \HH (B_{n+1})$. In the case of the singular Hecke algebras, the natural embedding $SB_n 
\hookrightarrow SB_{n+1}$ also leads to a homomorphism $\iota_n : \HH (SB_n) \to \HH (SB_{n+1})$, 
but we do not know whether this homomorphism is injective.

\bigskip\noindent
Let $z$ be a new variable. Let $d \ge 0$. A {\it Markov trace} on the sequence $\{ \HH 
(S_dB_n)\}_{n=1}^{+\infty}$ is defined to be a collection of $\K$-linear maps
\[
\tr_n^d: \HH (S_dB_n) \to \K (z)\,, \quad n \ge 1\,,
\]
such that
\begin{itemize}
\item
$\tr_n^d (\alpha \beta) = \tr_n^d (\beta \alpha)$ for all singular braids $\alpha \in S_kB_n$ and 
$\beta \in S_lB_n$ such that $k+l=d$, and all $n \ge 1$;
\item
$\tr_{n+1}^d \circ \iota_n = \tr_n^d$ for all $n \ge 1$;
\item
$\tr_{n+1}^d (\iota_n(\beta) \sigma_n) = z \cdot \tr_n^d (\beta)$ for all $\beta \in S_dB_n$ and all $n 
\ge 1$.
\end{itemize}
Define a {\it Markov trace} on the sequence $\{ \HH (SB_n) \}_{n=1}^{+\infty}$ to be a collection of 
$\K$-linear maps
\[
\tr_n: \HH (SB_n) \to \K (z)\,, \quad n \ge 1\,,
\]
such that
\begin{itemize}
\item
$\tr_n (\alpha \beta) = \tr_n (\beta \alpha)$ for all singular braids $\alpha ,\beta \in SB_n$, and all 
$n \ge 1$;
\item
$\tr_{n+1} \circ \iota_n = \tr_n$ for all $n \ge 1$;
\item
$\tr_{n+1} (\iota_n(\beta) \sigma_n) = z \cdot \tr_n (\beta)$ for all $\beta \in SB_n$ and all $n \ge 1$.
\end{itemize}

\bigskip\noindent
Note that, if $T= \{ \tr_n\}_{n=1}^{+\infty}$ is a Markov trace on $\{ \HH (SB_n) 
\}_{n=1}^{+\infty}$, then, for all $d \ge 0$, the collection  $T^d = \{ \tr_n^d = 
\tr_n|_{ \HH (S_dB_n)} \}_{n=1}^{+\infty}$ of restrictions is a Markov trace on $\{ \HH 
(S_dB_n)\}_{n=1}^{+\infty}$. Conversely, a collection $\{T^d\}_{d=0}^{+\infty}$, where $T^d$ is a 
Markov trace on $\{ \HH(S_dB_n)\}_{n=1}^{+\infty}$ for all $d \ge 0$, determines a unique Markov trace on $\{ \HH 
(SB_n)\}_{n=1}^{+\infty}$ . So, both definitions of Markov traces are more or less equivalent. Now, 
since the number $d$ of singular points can be fixed in our study, we will mainly consider Markov 
traces with a fixed number of singular points in the remainder.

\bigskip\noindent
{\bf Remark.} We do not impose the condition $\tr_1 (1)=1$ in the above definitions because this 
condition has no real meaning in the context of singular braids. Moreover, without this condition, the 
Markov traces on $\{ \HH(SB_n)\}_{n=1}^{+\infty}$ (or on $\{ \HH(S_dB_n)\}_{n=1}^{+\infty}$) form a 
$\K(z)$-vector space. This will be of importance in the remainder. 

\bigskip\noindent
Let $\LL_d$ denote the set of (isotopy classes of) singular links with $d$ singular points. We fix a 
Markov trace $T= \{ tr_n^d \}_{n=1}^{+\infty}$ on $\{ \HH(S_dB_n)\}_{n=1}^{+\infty}$, and turn to 
define an invariant $I_T : \LL_d \to \K (\sqrt{y})$. We follow the same strategy as Jones in 
\cite{Jones1}.

\bigskip\noindent
Let $\pi: SB_n \to \HH (SB_n)$ denote the natural map, and let $\varepsilon: SB_n \to \Z$ be the 
homomorphism defined by
\[
\varepsilon (\sigma_i)=1\,,\ \varepsilon(\sigma_i^{-1}) = -1\,,\ \varepsilon (\tau_i)=0\,,\ \text{for } 
1 \le i\le n-1\,.
\]
We consider the following change of variables:
\[
z=\frac{q-1}{1-qy} \quad \Leftrightarrow \quad y = \frac{z-q+1}{qz}\,.
\]
For a braid $\beta \in S_dB_n$ we set
\[
I_T(\beta) = \left( \frac{q-1}{1-qy} \right)^{-n+1} \cdot (\sqrt{y})^{\varepsilon(\beta) -n+1} \cdot \tr_n^d 
(\pi (\beta))\,.
\]
This is an element of $\K (\sqrt{y})$.

\bigskip\noindent
{\bf Proposition 3.3.} {\it Let $(\alpha,n)$ and $(\beta,m)$ be two singular braids with $d$ singular 
points. If $\hat \alpha$ is isotopic to $\hat \beta$, then $I_T(\alpha) = I_T(\beta)$.}

\bigskip\noindent
{\bf Proof.} By Theorem 2.4, it suffices to consider the following cases:
\begin{enumerate}
\item
$n=m$, and there exist $\gamma_1 \in S_kB_n$, $\gamma_2 \in S_lB_n$ such that $k+l=d$, $\alpha= 
\gamma_1 \gamma_2$, and $\beta= \gamma_2 \gamma_1$;
\item
$m=n+1$ and $\beta = \alpha \sigma_n$;
\item
$m=n+1$ and $\beta= \alpha \sigma_n^{-1}$.
\end{enumerate}

\bigskip\noindent
Suppose that $n=m$ and there exist $\gamma_1 \in S_kB_n$, $\gamma_2 \in S_lB_n$ such that $k+l=d$, 
$\alpha = \gamma_1 \gamma_2$, and $\beta= \gamma_2 \gamma_1$. Then, by definition, we have $\tr_n^d 
(\pi(\alpha)) = \tr_n^d (\pi(\beta))$ and $\varepsilon (\alpha) = \varepsilon(\beta)$, thus $I_T(\alpha) = I_T 
(\beta)$.

\bigskip\noindent
Suppose that $m=n+1$ and $\beta = \alpha \sigma_n$. Then 
\[\begin{array}{rcl}
I_T(\beta) &=& \left( \frac{q-1}{1-qy} \right)^{-m+1} \cdot (\sqrt{y})^{\varepsilon (\beta) -m +1} 
\cdot \tr_m^d(\pi(\beta))\\
&=& \left( \frac{q-1}{1-qy} \right)^{-n} \cdot (\sqrt{y})^{\varepsilon (\alpha) -n +1} \cdot 
\tr_{n+1}^d(\pi(\alpha)\sigma_n)\\
&=& \left( \frac{q-1}{1-qy} \right)^{-n} \cdot (\sqrt{y})^{\varepsilon (\alpha) -n +1} \cdot \left( 
\frac{q-1}{1-qy} \right) \cdot \tr_n^d(\pi(\alpha))\\
&=& I_T(\alpha)\,.
\end{array}\]

\bigskip\noindent
Suppose that $m=n+1$ and $\beta= \alpha \sigma_n^{-1}$. Recall the equality
\[
\sigma_n^{-1} = q^{-1} \sigma_n - q^{-1} (q-1)\,.
\]
Then
\[\begin{array}{rcl}
I_T(\beta) &=& \left( \frac{q-1}{1-qy} \right)^{-m+1} \cdot (\sqrt{y})^{\varepsilon (\beta) -m+1} \cdot 
\tr_m^d (\pi (\beta))\\
&=& \left( \frac{q-1}{1-qy} \right)^{-n} \cdot (\sqrt{y})^{\varepsilon (\alpha) -n-1} \cdot \tr_{n+1}^d 
(\pi (\alpha) \sigma_n^{-1})\\
&=& \left( \frac{q-1}{1-qy} \right)^{-n} \cdot (\sqrt{y})^{\varepsilon (\alpha) -n-1} \cdot \left( q^{-
1} \tr_{n+1}^d (\pi (\alpha) \sigma_n) - q^{-1}(q-1) \tr_{n+1}^d (\pi (\alpha)) \right)\\
&=& \left( \frac{q-1}{1-qy} \right)^{-n} \cdot (\sqrt{y})^{\varepsilon (\alpha) -n-1} \cdot \left( 
\frac{q-1}{1-qy} \right) \cdot y \cdot \tr_{n}^d (\pi (\alpha))\\
&=& I_T(\alpha)\,.
\end{array}\]
\qed

\bigskip\noindent
For $L \in \LL_d$, we choose a singular braid $(\beta,n)$ such that $\hat \beta = L$, and we set 
$I_T(L)= I_T(\beta)$. By Proposition 3.3, $I_T(L)$ is a well-defined invariant.

\bigskip\noindent
Let $A$ be an abelian group, let $I: \LL_d \to A$ be an invariant, and let $t,x \in A$. We say that $I$ 
satisfies the {\it $(t,x)$ skein relation} if
\[
t^{-1} \cdot I(L_+) -t \cdot I(L_-) = x \cdot I(L_0)\,,
\]
for all singular links $L_+, L_-, L_0 \in \LL_d$ that have the same link diagram except in the 
neighborhood of a crossing where they are like in Figure 3.1.

\begin{figure}[htb]
\bigskip
\centerline{
\setlength{\unitlength}{.5cm}
\begin{picture}(15,4)
\put(0,1){\includegraphics[width=7.5cm]{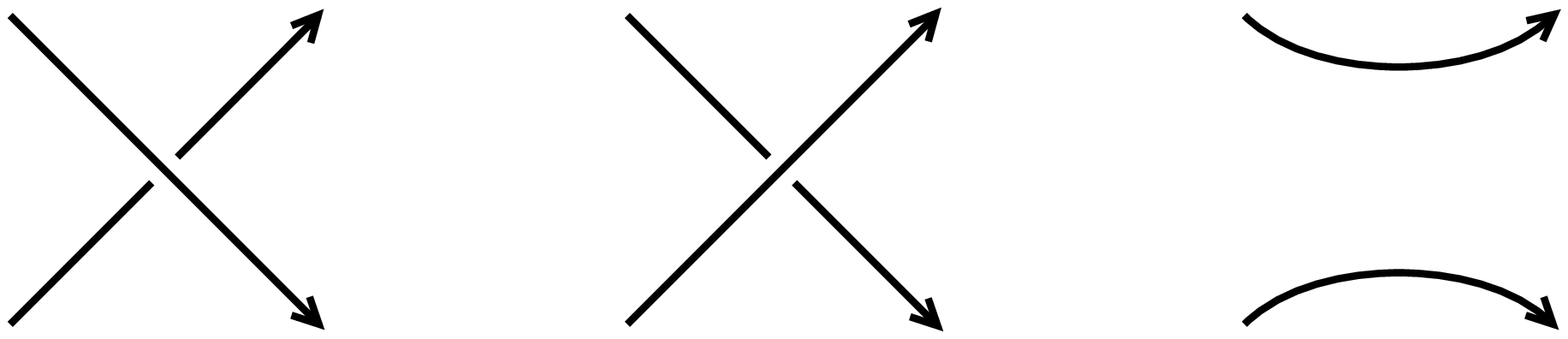}}
\put(1,0){\small $L_+$}
\put(7,0){\small $L_-$}
\put(13,0){\small $L_0$}
\end{picture}}
\bigskip
\centerline{{\bf Figure 3.1.}  The singular links $L_+$, $L_-$, and $L_0$.}
\end{figure}

\bigskip\noindent
Now, we set
\[
t = \sqrt{y} \cdot \sqrt{q}\,, \quad x= \sqrt{q} - \frac{1}{\sqrt{q}}\,,
\]
and we define $\widetilde{\tr}_n^d: S_dB_n \to \C (\sqrt{q}, \sqrt{y})$ by
\[
\widetilde{\tr}_n^d (\beta)= (\sqrt{q})^{-\varepsilon (\beta)} \cdot \tr_n^d (\pi (\beta))\,.
\]
With these new notations $I_T(\beta)$ can be written
\[
I_T(\beta)= \left( \frac{1-t^2}{tx} \right)^{n-1} \cdot t^{\varepsilon (\beta)} \cdot \widetilde{\tr}_n^d 
(\beta)\,.
\]

\bigskip\noindent
{\bf Proposition 3.4.} {\it The invariant $I_T: \LL_d \to \C(\sqrt{q},\sqrt{y})$ satisfies the 
$(t,x)$ skein relation.}

\bigskip\noindent
{\bf Proof.} 
Let $L_+, L_-, L_0 \in \LL_d$ be three singular links that have the same 
singular link diagram except in the neighborhood of a non-singular crossing where they are like in 
Figure 3.1.
A careful reading of the proof of Theorem 2.3 shows that there exist a singular braid 
$(\beta,n)$ with $d$ singular points, and an index $1 \le i\le n-1$, such that $L_+= \widehat{\beta 
\sigma_i}$, $L_- = \widehat{\beta \sigma_i^{-1}}$, and $L_0= \hat \beta$. On the other hand, the Hecke 
relation \eqref{R31} implies
\[
\widetilde{\tr}_n^d (\beta \sigma_i) = x \widetilde{\tr}_n^d (\beta) + \widetilde{\tr}_n^d (\beta 
\sigma_i^{-1})\,.
\]
Hence,
\[\begin{array}{rl}
&t^{-1} I_T(L_+) -t I_T(L_-)\\
=&\left( \frac{1-t^2}{tx}\right)^{n-1} t^{\varepsilon (\beta)} (\widetilde{\tr}_n^d (\beta \sigma_i) - 
\widetilde{\tr}_n^d (\beta \sigma_i^{-1}))\\
=&x\left( \frac{1-t^2}{tx}\right)^{n-1} t^{\varepsilon (\beta)} \widetilde{\tr}_n^d (\beta)\\
=& I_T (L_0)\,.
\end{array}\]
\qed

\bigskip\noindent
Define a {\it Markov trace} on the sequence $\{ \HH (S_dB_n) \}_{n=1}^{+\infty}$ {\it with coefficients 
in $\C (\sqrt{q}, \sqrt{y})$} to be a collection of $\K$-linear maps
\[
\tr_n^d: \HH (S_dB_n) \to \C (\sqrt{q}, \sqrt{y})\,, \quad n\ge 1\,,
\]
such that
\begin{itemize}
\item
$\tr_n^d (\alpha \beta) = \tr_n^d (\beta \alpha)$ for all singular braids $\alpha \in S_kB_n$ and 
$\beta \in S_lB_n$ such that $k+l=d$, and all $n \ge 1$;
\item
$\tr_{n+1}^d \circ \iota_n = \tr_n^d$ for all $n \ge 1$;
\item
$\tr_{n+1}^d (\iota_n(\beta) \sigma_n) = z \cdot \tr_n^d (\beta)$ for all $\beta \in S_dB_n$ and all $n 
\ge 1$.
\end{itemize}
Using the same trick as above, a Markov trace $T$ on the sequence $\{ \HH (S_dB_n)\}_{n=1}^{+\infty}$ 
with coefficients in $\C (\sqrt{q}, \sqrt{y})$
defines an invariant $I_T: \LL_d \to \C (\sqrt{q}, \sqrt{y})$ which satisfies the $(t,x)$ skein relation for 
$t= \sqrt{y} \sqrt{q}$ and $x= \sqrt{q} - \frac{1}{\sqrt{q}}$. Now, the reverse of Proposition 3.4 is 
true in the following sense.

\bigskip\noindent
{\bf Proposition 3.5.} {\it Let $I: \LL_d \to \C (\sqrt{q}, \sqrt{y})$ be an invariant which satisfies 
the $(t,x)$ skein relation for $t= \sqrt{y} \sqrt{q}$ and $x= \sqrt{q} - \frac{1}{\sqrt{q}}$. Then there exists 
a Markov trace $T$ on $\{ \HH (S_dB_n) \}_{n=1}^{+\infty}$ with coefficients in $\C( \sqrt{q}, 
\sqrt{y})$ such that $I=I_T$.}

\bigskip\noindent
{\bf Proof.} Recall that $\pi: SB_n \to \HH (SB_n)$ denotes the natural map, and that $\varepsilon: 
SB_n \to \Z$ is the homomorphism defined by $\varepsilon (\sigma_i) = 1$, $\varepsilon (\sigma_i^{-1}) 
= -1$, and $\varepsilon(\tau_i)=0$, for $1 \le i\le n-1$.

\bigskip\noindent
Let $\widetilde{\tr}_n^d : S_dB_n \to \C (\sqrt{q}, \sqrt{y})$ be the map defined by
\[
\widetilde{tr}_n^d (\beta)= t^{-\varepsilon (\beta)} \cdot \left( \frac{1-t^2}{tx} \right)^{1-n} \cdot 
I (\hat \beta)\,.
\]
Let $\alpha \in S_kB_n$ and $\beta \in S_lB_n$ such that $k+°l=d$, and let $1 \le i\le n-1$. By the $(t,x)$
skein relation we have
\[
t^{-1} \cdot I( \widehat{ \alpha \sigma_i^2 \beta} )- t \cdot I( \widehat{ \alpha \beta}) = x \cdot 
I( \widehat{ \alpha \sigma_i \beta})\,,
\]
thus
\begin{equation}\label{R36}
\widetilde{tr}_n^d (\alpha \sigma_i^2 \beta) = x \widetilde{tr}_n^d (\alpha \sigma_i \beta) + 
\widetilde{tr}_n^d( \alpha \beta)\,.
\end{equation}
Let $\alpha \in S_kB_n$ and $\beta \in S_lB_n$ such that $k+l=d$. We have $\widehat{ \alpha \beta} = 
\widehat{ \beta \alpha}$, thus
\begin{equation}\label{R37}
\widetilde{\tr}_n^d (\alpha \beta) = \widetilde{\tr}_n^d (\beta \alpha)\,.
\end{equation}
Let $\beta \in S_dB_n$. We have $\widehat{\beta \sigma_n} = \widehat{\beta \sigma_n^{-1}} = \hat 
\beta$, thus, by the $(t,x)$ skein relation,
\[\begin{array}{rc}
&t^{-1} \cdot I( \widehat{ (\beta \sigma_n, n+1)} )- t \cdot I( \widehat{ (\beta \sigma_n^{-1}, n+1)} 
)= x \cdot I( \widehat{ (\beta, n+1)} )\\
\Rightarrow\quad& t^{-1} (1-t^2) \cdot I( \widehat{ (\beta,n)}) = x\cdot I( \widehat{ (\beta,n+1)})\,,
\end{array}\]
therefore
\begin{equation}\label{R38}
\widetilde{\tr}_n^d (\beta) = \widetilde{\tr}_{n+1}^d (\beta)\,.
\end{equation}
Let $\beta \in S_dB_n$. Since $\widehat{\beta \sigma_n} = \hat \beta$, we have $I( \widehat{( \beta 
\sigma_n, n+1)}) = I( \widehat{( \beta,n)})$, thus
\begin{equation}\label{R39}
\widetilde{\tr}_{n+1}^d (\beta \sigma_n) = \left( \frac{x}{1-t^2} \right) \cdot \widetilde{\tr}_n^d 
(\beta)\,.
\end{equation}

\bigskip\noindent
Let $\tr_n^u: \K [S_dB_n] \to \C( \sqrt{q}, \sqrt{y})$ be the $\K$-linear map defined by
\[
\tr_n^u (\beta)= (\sqrt{q})^{\varepsilon (\beta)} \cdot \widetilde{\tr}_n^d (\beta)\,, \quad\text{for } 
\beta \in S_dB_n\,.
\]
Let $\alpha \in S_kB_n$ and $\beta \in S_lB_n$ such that $k+°l=d$, and let $1 \le i\le n-1$. By 
\eqref{R36}, we have
\[
\tr_n^u (\alpha \sigma_i^2 \beta) = (q-1) \cdot \tr_n^u (\alpha \sigma_i \beta) + q \cdot \tr_n^u 
(\alpha \beta)\,.
\]
So, $\tr_n^u : \K [S_dB_n] \to \C (\sqrt{q}, \sqrt{y})$ induces a $\K$-linear map $\tr_n^d: \HH 
(S_dB_n) \to \C (\sqrt{q}, \sqrt{y})$.

\bigskip\noindent
The relations \eqref{R37}, \eqref{R38}, and \eqref{R39} imply that $T= \{ \tr_n^d \}_{n=1}^{+\infty}$ 
is a Markov trace on $\{ \HH (S_dB_n)\}_{n=1}^{+\infty}$ with coefficients in $\C (\sqrt{q}, 
\sqrt{y})$, and a direct calculation shows that $I=I_T$.
\qed


\section{The space of traces}

For $d \ge 0$, we denote by $\TR_d$ the set of all traces on $\{ \HH (S_dB_n) \}_{n=1}^{+\infty}$. This 
is a $\K (z)$-vector space. Note also that the space of all traces on $\{\HH (SB_n) \}_{n=1}^{+\infty}$ 
is the completion of $\TR= \oplus_{d=0}^{+\infty} \TR_d$. We start our analysis recalling the 
following.

\bigskip\noindent
{\bf Theorem 4.1} (Ocneanu \cite{Jones1}, \cite{HOMFLY}). {\it There exists a unique trace $T_0^0 = \{ 
\tr_n^0\}_{n=1}^{+\infty}$ on $\{ \HH (B_n) \}_{n=1}^{+\infty}$ such that $\tr_1^0 (1)=1$.}
\qed

\bigskip\noindent
{\bf Corollary 4.2.} {\it $\TR_0$ is a one-dimensional $\K(z)$-vector space spanned by $T_0^0$.}
\qed

\bigskip\noindent
The above trace $T_0^0$ is called the {\it Ocneanu trace}. It will be a master piece in our study.

\bigskip\noindent
In this section we prove that $\TR_d$ is of dimension $d+1$ and construct an explicit 
basis $\{ T_0^d, T_1^d, \dots, 
\linebreak
T_d^d\}$ of $\TR_d$.

\bigskip\noindent
We start with the definition of the Markov traces $T_k^d$, $0 \le k \le d$.

\bigskip\noindent
Let
\[
g_0^d, g_1^d: \HH (S_{d+1}B_n) \to \HH (S_dB_n)
\]
be the $\K$-linear map defined as follows. Let $\beta \in S_{d+1}B_n$. Write $\beta$ in the form
\[
\beta= \alpha_0 \tau_{i_1} \alpha_1 \cdots \tau_{i_d} \alpha_d \tau_{i_{d+1}} \alpha_{d+1}\,,
\]
where $1 \le i_j \le n-1$ for $1 \le j \le d+1$, and $\alpha_j \in B_n$ for $0 \le j\le d+1$. Then
\begin{align*}
g_0^d (\beta)&= \sum_{j=1}^{d+1} \alpha_0 \tau_{i_1} \alpha_1 \cdots \tau_{i_{j-1}} \alpha_{j-1} \cdot 
\alpha_j \cdot \tau_{i_{j+1}} \alpha_{j+1} \cdots \tau_{i_{d+1}} \alpha_{d+1}\,,\\
g_1^d (\beta)&= \sum_{j=1}^{d+1} \alpha_0 \tau_{i_1} \alpha_1 \cdots \tau_{i_{j-1}} \alpha_{j-1} 
\cdot \sigma_{i_j} \alpha_j \cdot \tau_{i_{j+1}} \alpha_{j+1} \cdots \tau_{i_{d+1}} \alpha_{d+1}\,.
\end{align*}
It is easily seen from the presentation of $SB_n$ given in Theorem 2.2 that $g_0^d$ and $g_1^d$ are 
well-defined.

\bigskip\noindent
Let 
\[
\Phi_0^d, \Phi_1^d: \TR_d \to \TR_{d+1}
\]
be the $\K(z)$-linear maps defined as follows. Let $T= \{ \tr_n^d \}_{n=1}^{+\infty}$ be an element of 
$\TR_d$. Then, for $\omega \in \HH (S_{d+1}B_n)$, we set
\[
\Phi_0^d (T) (\omega) = \tr_n^d (g_0^d (w))\,, \quad \Phi_1^d (T) (\omega) = \tr_n^d (g_1^d (w))\,.
\]
It is easily checked that $\Phi_\varepsilon^d \circ \Phi_\mu^{d-1} = \Phi_\mu^d \circ 
\Phi_\varepsilon^{d-1}$ for all $\varepsilon, \mu \in \{0,1\}$, and all $d \ge 1$.

\bigskip\noindent
Now, we define $T_k^d$ by induction on $d$. According to the previous notation, $T_0^0$ is the Ocneanu 
trace of Theorem 4.1. If $d \ge 1$, then
\[
T_k^d= \left\{\begin{array}{ll}
\Phi_0^{d-1} (T_k^{d-1}) \quad&\text{if } k\le d-1\,,\\
\Phi_1^{d-1}(T_{d-1}^{d-1}) \quad &\text{if } k=d\,.
\end{array}\right.
\]
Note that we also have $T_k^d = \Phi_1^{d-1} (T_{k-1}^{d-1})$ for all $1 \le k\le d-1$.

\bigskip\noindent
{\bf Theorem 4.3.} {\it Let $d \ge 0$. Then $\{ T_0^d, T_1^d, \dots, T_d^d \}$ is a linearly 
independent family of $\TR_d$.}

\bigskip\noindent
The following lemmas 4.4 to 4.6 are preliminaries to the proof of Theorem 4.3.

\bigskip\noindent
The submonoid of $B_n$ generated (as a monoid) by $\sigma_1, \dots, \sigma_{n-1}$ is called the {\it 
positive braid monoid} and is denoted by $B_n^+$. By \cite{Garsi1}, it has a monoid presentation with 
generators $\sigma_1, \dots, \sigma_{n-1}$ and relations
\[\begin{array}{cl}
\sigma_i \sigma_j \sigma_i = \sigma_j \sigma_i \sigma_j &\quad\text{if } |i-j|=1\,,\\
\sigma_i \sigma_j = \sigma_j \sigma_i &\quad\text{if } |i-j| \ge 2\,.
\end{array}\]

\bigskip\noindent
{\bf Lemma 4.4.} {\it Let $n \ge 1$, and let $\beta \in B_n^+$. Then $T_0^0 (\beta) \in \Z [q,z]$.}

\bigskip\noindent
{\bf Proof.} Let $U_n$ denote the $\Z[q]$-submodule of $\HH (B_n)$ spanned by $B_n^+$. We prove by 
induction on $n \ge 2$ that $U_n$ is actually spanned as a $\Z[q]$-module by $B_{n-1}^+ \cup \{ \alpha 
\sigma_{n-1} \alpha'; \alpha, \alpha' \in B_{n-1}^+ \}$.

\bigskip\noindent
Suppose $n=2$. Then $U_2$ is spanned as a $\Z[q]$-module by $\{ \sigma_1^a; a \ge 0\}$. Now, the Hecke 
relation \eqref{R31} implies that
\[
\sigma_1^a = (q-1) \sigma_1^{a-1} + q \sigma_1^{a-2}\,, \quad\text{for all } a \ge 2\,,
\]
thus $U_2$ is spanned by $\{1, \sigma_1\}$.

\bigskip\noindent
Suppose $n  \ge 3$. Let $V_n$ be the $\Z[q]$-submodule spanned by $B_{n-1}^+ \cup \{ \alpha \sigma_{n-
1} \alpha'; \alpha, \alpha' \in B_{n-1}^+\}$. Let $\beta \in B_n^+$. We write $\beta$ in the form
\[
\beta= \beta_0 \sigma_{n-1} \beta_1 \cdots \sigma_{n-1} \beta_l\,,
\]
where $\beta_0, \beta_1, \dots, \beta_l \in B_{n-1}^+$, and prove that $\beta \in V_n$ by induction on 
$l$. The cases $l=0$ and $l=1$ are obvious. So, we can suppose that $l \ge 2$. By induction (on $n$), 
we can assume that either $\beta_1 \in B_{n-2}^+$, or $\beta_1 = \beta_1' \sigma_{n-2} \beta_1''$ for 
some $\beta_1', \beta_1'' \in B_{n-2}^+$. If $\beta_1 \in B_{n-2}^+$, then
\[\begin{array}{rl}
\beta &= \beta_0 \beta_1 \sigma_{n-1}^2 \beta_2 \sigma_{n-1} \beta_3 \cdots \sigma_{n-1} \beta_l \\
&= (q-1) \cdot \beta_0 \beta_1 \sigma_{n-1} \beta_2 \sigma_{n-1} \beta_3 \cdots \sigma_{n-1} \beta_l + 
q \cdot\beta_0 \beta_1 \beta_2 \sigma_{n-1} \beta_3 \cdots \sigma_{n-1} \beta_l\,,
\end{array}\]
thus, by induction (on $l$), we have $\beta \in V_n$. If $\beta_1 = \beta_1' \sigma_{n-2} \beta_1''$ 
for some $\beta_1', \beta_1'' \in B_{n-2}^+$, then
\[\begin{array}{rl}
\beta &= \beta_0 \beta_1' \sigma_{n-1} \sigma_{n-2} \sigma_{n-1} \beta_1'' \beta_2 \sigma_{n-1} \beta_3 
\cdots \sigma_{n-1} \beta_l\\ 
&= (\beta_0 \beta_1' \sigma_{n-2}) \sigma_{n-1} (\sigma_{n-2} \beta_1'' \beta_2) \sigma_{n-1} \beta_3 
\cdots \sigma_{n-1} \beta_l\,,
\end{array}\]
thus, by induction (on $l$), we have $\beta \in V_n$.

\bigskip\noindent
Now, we take $\beta \in B_n^+$ and turn to prove that $T_0^0( \beta) \in \Z [q,z]$. We argue by 
induction on $n$.

\bigskip\noindent
Suppose $n \ge 2$. By the above observation, we can assume that either $\beta \in B_{n-1}^+$, or $\beta 
= \alpha \sigma_{n-1} \alpha'$ for some $\alpha, \alpha' \in B_{n-1}^+$. If $\beta \in B_{n-1}^+$, then, 
by induction, $T_0^0 (\beta) \in \Z [q,z]$. If $\beta = \alpha \sigma_{n-1} \alpha'$ for some $\alpha, 
\alpha' \in B_{n-1}^+$, then, by induction, $T_0^0(\beta) = z \cdot T_0^0 (\alpha \alpha') \in \Z 
[q,z]$.
\qed

\bigskip\noindent
{\bf Lemma 4.5.} {\it Let $1 \le a \le n-1$, and let $\alpha, \alpha' \in \langle \sigma_{a+1}, \dots, 
\sigma_{n-1} \rangle^+$, where $\langle \sigma_{a+1}, \dots, \sigma_{n-1} \rangle^+$ denotes the 
submonoid generated by $\sigma_{a+1}, \dots, \sigma_{n-1}$. Then
\begin{gather*}
T_0^0 (\alpha \sigma_a \alpha')|_{z=0} = 0\,,\\
T_0^0 (\alpha \sigma_a^2 \alpha') |_{z=0} = q \cdot T_0^0 (\alpha \alpha') |_{z=0}\,.
\end{gather*}}

\bigskip\noindent
{\bf Proof.} The first equality is a consequence of the following one
\[
T_0^0( \alpha \sigma_a \alpha')= z \cdot T_0^0( \alpha \alpha')
\]
whose proof is left to the reader. The second equality follows from the first one and the Hecke 
relation \eqref{R31}.
\qed

\bigskip\noindent
The following lemma is a direct consequence of the previous one.

\bigskip\noindent
{\bf Lemma 4.6.} {\it Let $0 \le a,b \le n-1$, and let $i_1, \dots, i_a \in \{1, \dots, n-1 \}$ such 
that $i_1 <i_2< \cdots <i_a$. Let
\[
\gamma= \sigma_{i_a} \cdots \sigma_{i_2} \sigma_{i_1} \sigma_1 \sigma_2 \cdots \sigma_b\,.
\]
Then
\[
T_0^0(\gamma)|_{z=0}= \left\{\begin{array}{ll}
q^a &\quad\text{if } a=b,\, i_1=1,\, \dots,\, i_a=a\,,\\
0& \quad \text{otherwise}\,.
\end{array}\right.
\]}
\qed

\bigskip\noindent
{\bf Proof of Theorem 4.3.} For $0 \le b \le d$, we set
\[
\gamma_b^d= \tau_d \cdots \tau_2 \tau_1 \sigma_1 \sigma_2 \cdots \sigma_b\,.
\]
A direct calculation shows that
\[
T_a^d (\gamma_b^d)= (d-a)!\, a! \sum_{1 \le i_1< \cdots <i_a \le n-1} T_0^0 (\sigma_{i_a} \cdots 
\sigma_{i_2} \sigma_{i_1} \sigma_1 \sigma_2 \cdots \sigma_b)\,.
\]
By Lemma 4.4, we have $T_a^d (\gamma_b^d) \in \Z [q,z]$, and, by Lemma 4.6,
\[
T_a^d (\gamma_b^d)|_{z=0} = \left\{
\begin{array}{ll}
(d-a)!\, a!\, q^a &\quad \text{if } a=b\,,\\
0 &\quad \text{otherwise}\,.
\end{array}\right.
\]
This implies that $T_0^d, T_1^d, \dots, T_d^d$ are linearly independent.
\qed

\bigskip\noindent
{\bf Theorem 4.7.} {\it Let $d \ge 0$. Then $\TR_d$ is a $\K(z)$-vector space of dimension $d+1$. In 
particular, $\{ T_0^d, T_1^d, \dots, T_d^d\}$ is a basis of $\TR_d$.}

\bigskip\noindent
The main ingredient in the proof of Theorem 4.7 are the relations in $\HH (SB_n)$ that will be proved in 
the following lemmas 4.8 to 4.11. We will prove Theorem 4.7 after these lemmas.

\bigskip\noindent
{\bf Lemma 4.8.} {\it Let $i,j \in \{1, \dots, n-1\}$ such that $|i-j|=1$, and let $a \ge 1$. Set
\[
B_{i\,j}= \sigma_i + \sigma_j -(q-1)\,.
\]
Then
\begin{gather}
\label{R41}
\sigma_i \tau_j^a = q^{-1} \tau_j^a \sigma_i \sigma_j \sigma_i - q^{-1} \tau_i^a \sigma_i 
\sigma_j \sigma_i + \tau_j^a \sigma_i \,;\\
\label{R42}
\sigma_j \tau_i^a = \tau_j^a (\sigma_i + \sigma_j -(q-1)) - \tau_i^a (\sigma_i - (q-1))\,;\\
\label{R43}
B_{i\,j} \tau_i^a = \tau_j^a B_{i\,j}\,;\\
\label{R44}
\begin{array}{rl}
&\tau_i^a (\sigma_i \sigma_j + \sigma_j \sigma_i -(q-1) \sigma_i - (q-1) \sigma_j + (q^2-q+1))\\ 
= &\tau_j^a (\sigma_i \sigma_j + \sigma_j \sigma_i -(q-1) \sigma_i - (q-1) \sigma_j + (q^2-q+1))\,.
\end{array}
\end{gather}}

\bigskip\noindent
{\bf Proof.} The equalities \eqref{R41}, \eqref{R42}, and \eqref{R44} are proved in Section 3 (see 
\eqref{R32}, \eqref{R35}, \eqref{R34}). Since $\tau_i^a$ commutes with $\sigma_i -(q-1)$, the equality 
\eqref{R43} is a direct consequence of \eqref{R42}.
\qed

\bigskip\noindent
{\bf Lemma 4.9.} {\it Let $i,j \in \{1, \dots, n-1\}$ such that $|i-j|=1$. Then $B_{i\,j}$ is invertible 
in $\HH (B_n)$.}

\bigskip\noindent
{\bf Proof.} A direct calculation shows that
\[
-q^{-1} (q+1)^{-2} (q(q-1) - 2q \sigma_i - 2q \sigma_j -(q-1) \sigma_i \sigma_j 
-(q-1) \sigma_j \sigma_i + 2\sigma_i \sigma_j \sigma_i)
\]
is the inverse of $B_{i\,j}$.
\qed

\bigskip\noindent
{\bf Lemma 4.10.} {\it Let $i,j \in \{1, \dots, n-1\}$ such that $|i-j| \ge 2$, and let $a \ge 1$. Let
\[
C_{i\,j}= 2 \sigma_i \sigma_j -(q-1) \sigma_i -(q-1)\sigma_j +q^2+1\,.
\]
Then
\begin{gather}
\label{R45}
\sigma_i \tau_j^a = \tau_j^a \sigma_i\,;\\
\label{R46}
(\sigma_i-\sigma_j)^2 = (q+1)^2 - C_{i\,j}\,;\\
\label{R47}
\tau_i^a C_{i\,j} = \tau_j^a C_{i\,j}\,. 
\end{gather}}

\bigskip\noindent
{\bf Proof.} The equality \eqref{R45} is a straightforward consequence of Theorem 2.2, and \eqref{R46} 
can be easily proved with a direct calculation. So, it remains to prove \eqref{R47}.

\bigskip\noindent
First, we study the case where $i=1$ and $j=3$. We apply twice \eqref{R42} to $\sigma_1 \sigma_3 
\tau_2^a$ and obtain
\[\begin{array}{rcl}
\sigma_1 \sigma_3 \tau_2^a &=& \sigma_1 \tau_3^a (\sigma_2 + \sigma_3 - (q-1)) - \sigma_1 \tau_2^a 
(\sigma_2 - (q-1))\\
&=& \tau_3^a \sigma_1 (\sigma_2 + \sigma_3 -(q-1)) - \tau_1^a (\sigma_1 + \sigma_2 -(q-1)) (\sigma_2 -
(q-1)) + \tau_2^a (\sigma_2 -(q-1))^2\\
&=& \tau_3^a (\sigma_1 \sigma_2 + \sigma_1 \sigma_3 - (q-1) \sigma_1) - \tau_1^a (\sigma_1 \sigma_2 - 
(q-1) \sigma_1 - (q-1) \sigma_2 + (q^2 -q+1))\\ 
&&\quad  + \tau_2^a (\sigma_2 -(q-1))^2\,.
\end{array}\]
Similarly,
\begin{multline*}
\sigma_3 \sigma_1 \tau_2^a = \tau_1^a (\sigma_3 \sigma_2 + \sigma_1 \sigma_3 - (q-1) \sigma_3) \\
- \tau_3^a (\sigma_3 \sigma_2 - (q-1) \sigma_2 - (q-1) \sigma_3 + (q^2 -q+1)) + \tau_2^a (\sigma_2 -(q-
1))^2\,.
\end{multline*}
Since $\sigma_1 \sigma_3 \tau_2^a = \sigma_3 \sigma_1 \tau_2^a$, it follows that
\begin{equation}\label{R48}
\begin{array}{rl}
&\tau_1^a (\sigma_1 \sigma_2 + \sigma_3 \sigma_2 + \sigma_1 \sigma_3 - (q-1) \sigma_1 - (q-1) \sigma_2 
- (q-1) \sigma_3 + (q^2-q+1))\\
=&\tau_3^a (\sigma_1 \sigma_2 + \sigma_3 \sigma_2 + \sigma_1 \sigma_3 - (q-1) \sigma_1 - (q-1) \sigma_2 
- (q-1) \sigma_3 + (q^2-q+1))\,.
\end{array}
\end{equation}
Set
\[
\omega_0 = \sigma_1 \sigma_2 + \sigma_3 \sigma_2 + \sigma_1 \sigma_3 - (q-1) \sigma_1 - (q-1) \sigma_2 - 
(q-1) \sigma_3 + (q^2-q+1)\,.
\]
By \eqref{R48}, we have $\tau_1^a \omega_0 = \tau_3^a \omega_0$. A direct calculation shows that
\[
C_{1\,3} = q^{-2} (q \sigma_1 \omega_0 + q \omega_0 \sigma_3 + (q-1) \sigma_1 \omega_0 \sigma_3 - 
\sigma_1 \omega_0 \sigma_3 \sigma_2)( \sigma_1 - (q-1))\,.
\]
Since $\tau_1$ and $\tau_3$ commute with $\sigma_1$, it follows that
\begin{equation}\label{R49}
\tau_1^a C_{1\,3} = \tau_3^a C_{1\,3}\,.
\end{equation}

\bigskip\noindent
Now, suppose that $1 \le i <j-1 \le n-2$. Set
\[
B_{j\,3} = B_{j\,j-1} \cdots B_{5\,4} B_{4\,3}\,, \quad B_{i\,1} = B_{i\,i-1} \cdots B_{3\,2} B_{2\,1}\,, 
\quad \delta_{i\,j} = B_{i\,1} B_{j\,3}\,.
\]
By Lemma 4.9, $\delta_{i\,j}$ is invertible, and by \eqref{R43}, we have
\[
\delta_{i\,j} \sigma_1 \delta_{i\,j}^{-1} = \sigma_i\,, \quad \delta_{i\,j} \tau_1 \delta_{i\,j}^{-1} = 
\tau_i\,, \quad \delta_{i\,j} 
\sigma_3 \delta_{i\,j}^{-1} = \sigma_j\,, \quad \delta_{i\,j} \tau_3 \delta_{i\,j}^{-1} = \tau_j\,,
\]
thus, by \eqref{R49},
\[
\tau_i^a C_{i\,j} = \delta_{i\,j} \tau_1^a C_{1\,3} \delta_{i\,j}^{-1} = \delta_{i\,j} \tau_3^a 
C_{1\,3} \delta_{i\,j}^{-1} = 
\tau_j^a C_{i\,j}\,.
\]
\qed

\bigskip\noindent
{\bf Lemma 4.11.} {\it Let $a,b \ge 1$. Then
\begin{equation}\label{R410}
\tau_1^a \tau_3^b (\sigma_3 - \sigma_1) = (\tau_2^b \tau_1^a + \tau_2^a \tau_3^b) (\sigma_3 - \sigma_1) 
+ \tau_2^{a+b} (B_{1\,2} - B_{2\,3})\,.
\end{equation}}

\bigskip\noindent
{\bf Proof.} Applying twice \eqref{R42} to $\sigma_2 \tau_1^a \tau_3^b$ we obtain
\[\begin{array}{rcl}
\sigma_2 \tau_1^a \tau_3^b &=& \tau_2^a (\sigma_1 +\sigma_2 -(q-1)) \tau_3^b - \tau_1^a (\sigma_1 -(q-
1)) \tau_3^b\\
&=&\tau_2^a \tau_3^b (\sigma_1 -(q-1)) + \tau_2^a \sigma_2 \tau_3^b - \tau_1^a \tau_3^b (\sigma_1 -(q-
1))\\
&=& \tau_2^a \tau_3^b (\sigma_1 -(q-1)) + \tau_2^{a+b} (\sigma_2 + \sigma_3 - (q-1)) - \tau_2^a 
\tau_3^b (\sigma_3 -(q-1))\\ 
&&\quad - \tau_1^a \tau_3^b (\sigma_1 -(q-1))\\
&=& \tau_2^a \tau_3^b (\sigma_1 -\sigma_3) + \tau_2^{a+b} B_{2\,3} - \tau_1^a \tau_3^b (\sigma_1 -(q-
1))\,.
\end{array}\]
Similarly,
\[
\sigma_2 \tau_3^b \tau_1^a = \tau_2^b \tau_1^a (\sigma_3 -\sigma_1) + \tau_2^{a+b} B_{1\,2} - \tau_1^a 
\tau_3^b (\sigma_3 -(q-1))\,.
\]
Since $\sigma_2 \tau_1^a \tau_3^b = \sigma_2 \tau_3^b \tau_1^a$, it follows that
\[
\tau_1^a \tau_3^b (\sigma_3 - \sigma_1) = (\tau_2^b \tau_1^a + \tau_2^a \tau_3^b) (\sigma_3 - \sigma_1) 
+ \tau_2^{a+b} (B_{1\,2} - B_{2\,3})\,.
\]
\qed

\bigskip\noindent
{\bf Proof of Theorem 4.7.} We fix once for all the number $d \ge 1$ of singular points. We set
\[
\widetilde{\TR}_d = \bigoplus_{n=2}^{+\infty}\big( \K(z) \otimes \HH(S_dB_n)\big)\,.
\]
Then $\TR_d$ can and will be viewed as the quotient of $\widetilde{TR}_d$ by the following relations:
\begin{itemize}
\item
$(\alpha \beta,n) = (\beta \alpha,n)$ for all $\alpha \in S_kB_n$ and $\beta \in S_lB_n$ such that 
$k+l=d$, and all $n\ge 2$;
\item
$(\beta,n+1)=(\beta,n)$ for all $\beta \in S_dB_n$ and all $n \ge 2$;
\item
$(\beta \sigma_n, n+1)= z \cdot (\beta,n)$ for all $\beta \in S_dB_n$ and all $n \ge 2$.
\end{itemize}
For $\omega \in \K (z) \otimes \HH(S_dB_n)$, we will denote by $[\omega]$ the element of $\TR_d$ 
represented by $\omega$.

\bigskip\noindent
We already know that $\dim TR_d \ge d+1$ (see Theorem 4.3). So, in order to prove Theorem 4.7, it 
suffices to show that $\TR_d$ is spanned by $d+1$ elements.

\bigskip\noindent
Recall the basis $\BB_n$ of $\HH(B_n)$ described in Section 3. For $n \ge 2$ we set
\[
\UU_n = \{ 1, \sigma_{n-1}, \sigma_{n-1} \sigma_{n-2}, \dots, \sigma_{n-1} \cdots \sigma_2 \sigma_1 
\}\,.
\]
Then $\BB_n$ is defined by induction on $n$ by
\[
\BB_1 = \{1\}\,, \quad \BB_n = \{ \beta u; \beta \in \BB_{n-1} \text{ and } u \in \UU_n \} \quad\text{if } n 
\ge 2\,.
\]
Let $\CC_n$ be the set of elements of $\TR_d$ of the form $[ \tau_{i_1} \cdots \tau_{i_d} \beta]$, where 
$1 \le i_j \le n-1$ for $1 \le j\le d$, and $\beta \in \BB_n$. Set $\CC_\infty = \cup_{n=2}^{+\infty} 
\CC_n$. By Proposition 3.1, $\CC_\infty$ spans $TR_d$.

\bigskip\noindent
Let $\omega \in \CC_n$. 
Let $1 \le l \le d$. If $\omega$ can be written in the form $\omega= [ \tau_{i_1}^{a_1} \cdots
\tau_{i_l}^{a_l} \beta]$, where $1 \le i_j \le n-1$ and $a_j \ge 1$ for $1 \le j\le l$, $a_1 + \cdots +
a_l=d$, and $\beta \in \BB_n$, then we say that {\it $\omega$ has a syllable length less or equal to $l$},
and we write $\Syl (\omega) \le l$. We set
\[
\DD_{l,n}= \{ \omega \in \CC_n; \Syl(\omega) \le l\}\,,\quad \text{and }\DD_{l,\infty} = \cup_{n=2}^{+\infty} 
\DD_{l,n}\,.
\]
Note that $\DD_{d,\infty} = \CC_\infty$ spans $\TR_d$.

\bigskip\noindent
For $\XX \subset \TR_d$, we denote by $\Span (\XX)$ the $\K(z)$-linear subspace spanned by $\XX$. The 
first step in the proof of Theorem 4.7 will consist on proving that $\Span (\DD_{l,\infty}) = 
\Span(\DD_{l-1,\infty})$ for all $l \ge 3$ (see Claims 1 to 4). Since $\TR_d= \Span (\DD_{d,\infty})$, 
it will follow that $\TR_d = \Span (\DD_{2, \infty})$. The second step will consist on proving that 
there exists a subset $\FF_3 \subset \DD_{2,3}$ with $d+1$ elements such that $\Span (\FF_3) = \Span 
(\DD_{2,\infty}) = \TR_d$ (see Claims 5 to 7).

\bigskip\noindent
Let $1 \le l \le d$, and let $1 \le r \le l$. Set $\varepsilon=2$ if $r$ is even, and $\varepsilon=1$ 
if $r$ is odd. Then we denote by $\EE_{r,l,n}$ the set of elements of $\DD_{l,n}$ of the form
\[
\omega= [ \tau_1^{a_1} \tau_2^{a_2} \tau_1^{a_3} \cdots \tau_\varepsilon^{a_r} \tau_{i_{r+1}}^{a_{r+1}} 
\cdots \tau_{i_l}^{a_l} \beta]\,,
\]
where $1 \le i_j \le n-1$ for $r+1 \le j \le l$, and $\beta \in \BB_n$. We set 
$\EE_{r,l,\infty} = \cup_{n=2}^{+\infty} \EE_{r,l,n}$.

\bigskip\noindent
{\bf Claim 1.} {\it Let $2 \le l \le d$. Then
\begin{equation}\label{R411}
\Span (\DD_{l,\infty}) = \Span (\DD_{l-1,\infty} \cup \EE_{2,l,\infty})\,.
\end{equation}}

\bigskip\noindent
{\bf Proof.} For $k \ge 2$, we denote by $\EE_{1,l,n}' (k)$ the set of elements $\omega \in 
\DD_{l,n}$ of the form $\omega= [\tau_{i_1}^{a_1} \tau_{i_2}^{a_2} \cdots \tau_{i_l}^{a_l} \beta]$, 
where $1 \le i_1\le k$, $1 \le i_j \le n-1$ for $2 \le j\le l$, and $\beta \in 
\BB_n$. We set $\EE_{1,l,\infty}'(k) = \cup_{n=2}^{+\infty} \EE_{1,l,n}' (k)$. Note that 
$\EE_{1,l,\infty}'(1) = \EE_{1,l,\infty}$, and $\DD_{l,n} = \DD_{l-1,n} \cup \EE_{1,l,n}' (n-1)$ for all 
$n \ge 2$.

\bigskip\noindent
We prove that
\begin{equation}\label{R412}
\Span (\DD_{l-1,\infty} \cup \EE_{1,l,\infty}' (k)) = \Span (\DD_{l-1,\infty} \cup \EE_{1,l,\infty}' 
(k-1))
\end{equation}
for all $k \ge 2$. This implies that
\begin{equation}\label{R413}
\Span( \DD_{l,\infty}) = \Span (\DD_{l-1,\infty} \cup \EE_{1,l,\infty})\,.
\end{equation}
Let $\omega \in \EE_{1,l,\infty}'(k)$ be of the form $\omega= [ \tau_k^{a_1} \tau_{i_2}^{a_2} 
\cdots \tau_{i_l}^{a_l} \beta]$. By \eqref{R41}, we have
\[\begin{array}{rcl}
\omega &=& [\tau_k^{a_1} \sigma_k \sigma_{k-1} \sigma_k \sigma_k^{-1} \sigma_{k-1}^{-1} 
\sigma_k^{-1} \tau_{i_2}^{a_2} \cdots \tau_{i_l}^{a_l} \beta]\\
&=& [\tau_{k-1}^{a_1} \tau_{i_2}^{a_2} \cdots \tau_{i_l}^{a_l} \beta] - q [\sigma_k\tau_{k-1}^{a_1} 
\sigma_k^{-1} \sigma_{k-1}^{-1} \sigma_k^{-1} \tau_{i_2}^{a_2} \cdots \tau_{i_l}^{a_l} \beta] + 
q [\tau_{k-1}^{a_1} \sigma_{k-1}^{-1} \sigma_k^{-1} \tau_{i_2}^{a_2} \cdots \tau_{i_l}^{a_l} 
\beta]\\
&=& [\tau_{k-1}^{a_1} \tau_{i_2}^{a_2} \cdots \tau_{i_l}^{a_l} \beta] - q [\tau_{k-1}^{a_1} 
\sigma_k^{-1} \sigma_{k-1}^{-1} \sigma_k^{-1} \tau_{i_2}^{a_2} \cdots \tau_{i_l}^{a_l} \beta 
\sigma_k]+ q [\tau_{k-1}^{a_1} \sigma_{k-1}^{-1} \sigma_k^{-1} \tau_{i_2}^{a_2} \cdots \tau_{i_l}^{a_l} 
\beta]\,.\\
\end{array}\]
It is easily checked by means of \eqref{R42} and \eqref{R45} that this element belongs to $\Span 
(\DD_{l-1,\infty} \cup \EE_{1,l,\infty}'(k-1))$.

\bigskip\noindent
Now, for $k\ge 2$, we denote by $\EE_{2,l,n}' (k)$ the set of elements $\omega \in \DD_{l,n}$ of 
the form $\omega=$ 
\linebreak
$[ \tau_1^{a_1} \tau_{i_2}^{a_2} \tau_{i_3}^{a_3} \cdots \tau_{i_l}^{a_l} \beta]$, 
where $2 \le i_2 \le k$, $1 \le i_j \le n-1$ for $3 \le j\le l$, and $\beta \in 
\BB_n$. We set $\EE_{2,l,\infty}' (k) = \cup_{n=2}^{+\infty} \EE_{2,l,n}'(k)$. Note that 
$\EE_{2,l,\infty}' (2)= \EE_{2,l,\infty}$, and $\EE_{2,l,n}'(n-1)= \EE_{1,l,n}$ for all $n \ge 2$.

\bigskip\noindent
Using the same arguments as in the proof of \eqref{R412}, one can easily show that
\begin{equation}\label{R414}
\Span (\DD_{l-1,\infty} \cup \EE_{2,l,\infty}' (k)) = \Span (\DD_{l-1, \infty} \cup \EE_{2,l,\infty}' 
(k-1))
\end{equation}
for all $k \ge 3$. (Here we also need to use the fact that $\sigma_k$ commutes with $\tau_1$.) It 
follows that
\[
\Span (\DD_{l,\infty}) = \Span (\DD_{l-1,\infty} \cup \EE_{2,l,\infty})\,.
\]
\qed

\bigskip\noindent
{\bf Claim 2.} {\it Let $l \ge 3$, and let $2 \le r\le l-1$. Then
\begin{equation}\label{R415}
\Span(\DD_{l-1,\infty} \cup \EE_{r,l,\infty}) = \Span( \DD_{l-1,\infty} \cup \EE_{r+1,l,\infty})\,.
\end{equation}}

\bigskip\noindent
{\bf Proof.} Set $\varepsilon=2$ if $r$ is even, and $\varepsilon=1$ if $r$ is odd. For $k\ge 3$, 
we denote by $\EE_{r+1,l,n}'(k)$ the set of elements $\omega \in \DD_{l,n}$ of the form $\omega= [ 
\tau_1^{a_1} \tau_2^{a_2} \tau_1^{a_3} \cdots \tau_\varepsilon^{a_r} \tau_{i_{r+1}}^{a_{r+1}} 
\tau_{i_{r+2}}^{a_{r+2}} \cdots \tau_{i_l}^{a_l} \beta]$, where $1 \le i_{r+1} \le 
k$, $1 \le i_j \le n-1$ for $r+2 \le j \le l$, and $\beta \in \BB_n$. We set $\EE_{r+1,l,\infty}'(k)= 
\cup_{n=2}^{+\infty} \EE_{r+1,l,n}' (k)$. Note that $\EE_{r+1,l,n}' (n-1)= \EE_{r,l,n}$ for all $n \ge 
2$.

\bigskip\noindent
Using the same arguments as in the proof of \eqref{R412}, one can easily show that 
\begin{equation}\label{R416}
\Span( \DD_{l-1,\infty} \cup \EE_{r+1,l,\infty}' (k)) = \Span( \DD_{l-1,\infty} \cup \EE_{r+1,l,\infty}' 
(k-1))
\end{equation}
for all $k \ge 4$. (Here we also need to use the fact that $\sigma_k$ commutes with $\tau_1$ and 
$\tau_2$.) This implies that
\begin{equation}\label{R417}
\Span (\DD_{l-1,\infty} \cup \EE_{r,l,\infty}) = \Span (\DD_{l-1,\infty} \cup \EE_{r+1,l,\infty}' 
(3))\,.
\end{equation}
Let $\omega \in \EE_{r+1,l,\infty}' (3)$ be an element of the form $\omega = [ \tau_1^{a_1} 
\tau_2^{a_2} \cdots \tau_\varepsilon^{a_r} \tau_3^{a_{r+1}} \tau_{i_{r+2}}^{a_{r+2}} \cdots 
\tau_{i_l}^{a_l} \beta]$. Now, in order to prove Claim 2, it suffices to show that such an element 
belongs to $\Span (\DD_{l-1,\infty} \cup \EE_{r+1,l,\infty})$.

\bigskip\noindent
Assume that $r$ is odd. So,
\[
\omega = [ \tau_1^{a_1} \cdots \tau_2^{a_{r-1}} \tau_1^{a_r} \tau_3^{a_{r+1}} \tau_{i_{r+2}}^{a_{r+2}} 
\cdots \tau_{i_l}^{a_l} \beta ]\,.
\]
Let
\[\begin{array}{rcl}
\omega_1&=& [ \tau_1^{a_1} \cdots \tau_2^{a_{r-1}} \tau_1^{a_r} \tau_3^{a_{r+1}} (\sigma_3 -\sigma_1)^2 
\tau_{i_{r+2}}^{a_{r+2}} \cdots \tau_{i_l}^{a_l} \beta ]\,,\\
\omega_2&=& [ \tau_1^{a_1} \cdots \tau_2^{a_{r-1}} \tau_1^{a_r} \tau_3^{a_{r+1}} C_{1\,3} 
\tau_{i_{r+2}}^{a_{r+2}} \cdots \tau_{i_l}^{a_l} \beta ]\,.
\end{array}\]
By Lemma 4.11, we have
\[\begin{array}{rcl}
\omega_1&=& [ \tau_1^{a_1} \cdots \tau_2^{a_{r-1}+a_r} \tau_3^{a_{r+1}} (\sigma_3 -\sigma_1)^2 
\tau_{i_{r+2}}^{a_{r+2}} \cdots \tau_{i_l}^{a_l} \beta ]\\ 
&& \quad + [ \tau_1^{a_1} \cdots \tau_2^{a_{r-1} +a_{r+1}} \tau_1^{a_r}  (\sigma_3 -\sigma_1)^2 
\tau_{i_{r+2}}^{a_{r+2}} \cdots \tau_{i_l}^{a_l} \beta ]\\
&&\quad + [ \tau_1^{a_1} \cdots \tau_2^{a_{r-1}+ a_r + a_{r+1}} (B_{1\,2} - B_{2\,3}) (\sigma_3 -
\sigma_1) \tau_{i_{r+2}}^{a_{r+2}} \cdots \tau_{i_l}^{a_l} \beta ]\\
&\in& \Span( \DD_{l-1, \infty}) \subset \Span (\DD_{l-1, \infty} \cup \EE_{r+1,l,\infty})\,.
\end{array}\]
On the other hand, by Lemma 4.10,
\[\begin{array}{rcl}
\omega_2&=& [ \tau_1^{a_1} \cdots \tau_2^{a_{r-1}} \tau_1^{a_r+a_{r+1}} C_{1\,3} 
\tau_{i_{r+2}}^{a_{r+2}} \cdots \tau_{i_l}^{a_l} \beta ]\\
&\in& \Span( \DD_{l-1, \infty}) \subset \Span (\DD_{l-1, \infty} \cup \EE_{r+1,l,\infty})\,.
\end{array}\]
Hence, by Lemma 4.10,
\[
\omega = (q+1)^{-2} (\omega_1+\omega_2) \in \Span (\DD_{l-1, \infty} \cup \EE_{r+1,l,\infty})\,.
\]

\bigskip\noindent
Now, assume that $r$ is even. So
\[
\omega = [ \tau_1^{a_1} \cdots \tau_1^{a_{r-1}} \tau_2^{a_r} \tau_3^{a_{r+1}} \tau_{i_{r+2}}^{a_{r+2}} 
\cdots \tau_{i_l}^{a_l} \beta]\,.
\]
Let
\[\begin{array}{rcl}
\omega_1 &=& [ \tau_1^{a_1} \cdots \tau_1^{a_{r-1}} \tau_2^{a_r} \tau_1^{a_{r+1}}B_{1\,2} B_{2\,3} 
\tau_{i_{r+2}}^{a_{r+2}} \cdots \tau_{i_l}^{a_l} \beta B_{1\,2}^{-2}]\,,\\
\omega_2 &=& [ \tau_2^{a_1} \cdots \tau_2^{a_{r-1}} \tau_1^{a_r} \tau_3^{a_{r+1}} (\sigma_1 - \sigma_3) 
\tau_{i_{r+2}}^{a_{r+2}} \cdots \tau_{i_l}^{a_l} \beta B_{1\,2}^{-1}]\,.
\end{array}\]
Obviously, $\omega_1 \in \Span (\DD_{l-1,\infty} \cup \EE_{r+1,l,\infty})$. On the other hand, by Lemma 
4.11,
\[\begin{array}{rcl}
\omega_2 &=& [ \tau_2^{a_1} \cdots \tau_2^{a_{r-1}+a_r} \tau_3^{a_{r+1}} (\sigma_1 - \sigma_3) 
\tau_{i_{r+2}}^{a_{r+2}} \cdots \tau_{i_l}^{a_l} \beta B_{1\,2}^{-1}]\\
&&\quad + [ \tau_2^{a_1} \cdots \tau_2^{a_{r-1}+a_{r+1}} \tau_1^{a_r} (\sigma_1 - \sigma_3) 
\tau_{i_{r+2}}^{a_{r+2}} \cdots \tau_{i_l}^{a_l} \beta B_{1\,2}^{-1}]\\
&&\quad + [ \tau_2^{a_1} \cdots \tau_2^{a_{r-1}+a_r +a_{r+1}} (B_{2\,3} - B_{1\,2}) 
\tau_{i_{r+2}}^{a_{r+2}} \cdots \tau_{i_l}^{a_l} \beta B_{1\,2}^{-1}]\\
&\in& \Span( \DD_{l-1, \infty}) \subset \Span (\DD_{l-1, \infty} \cup \EE_{r+1,l,\infty})\,.
\end{array}\]
Hence, by Lemma 4.8,
\[\begin{array}{rcl}
\omega &=& [ B_{1\,2} \tau_1^{a_1} \cdots \tau_1^{a_{r-1}} \tau_2^{a_r} \tau_3^{a_{r+1}} 
\tau_{i_{r+2}}^{a_{r+2}} \cdots \tau_{i_l}^{a_l} \beta B_{1\,2}^{-1}]\\
&=& [ \tau_2^{a_1} \cdots \tau_2^{a_{r-1}} \tau_1^{a_r} B_{1\,2} \tau_3^{a_{r+1}} 
\tau_{i_{r+2}}^{a_{r+2}} \cdots \tau_{i_l}^{a_l} \beta B_{1\,2}^{-1}]\\
&=& [ \tau_2^{a_1} \cdots \tau_2^{a_{r-1}} \tau_1^{a_r}  \tau_3^{a_{r+1}} (\sigma_1-(q-1)) 
\tau_{i_{r+2}}^{a_{r+2}} \cdots \tau_{i_l}^{a_l} \beta B_{1\,2}^{-1}]\\
&&\quad + [ \tau_2^{a_1} \cdots \tau_2^{a_{r-1}} \tau_1^{a_r} \sigma_2 \tau_3^{a_{r+1}} 
\tau_{i_{r+2}}^{a_{r+2}} \cdots \tau_{i_l}^{a_l} \beta B_{1\,2}^{-1}]\\
&=& [ \tau_2^{a_1} \cdots \tau_2^{a_{r-1}} \tau_1^{a_r}  \tau_3^{a_{r+1}} (\sigma_1-(q-1)) 
\tau_{i_{r+2}}^{a_{r+2}} \cdots \tau_{i_l}^{a_l} \beta B_{1\,2}^{-1}]\\
&&\quad + [ \tau_2^{a_1} \cdots \tau_2^{a_{r-1}} \tau_1^{a_r}\tau_2^{a_{r+1}} B_{2\,3}  
\tau_{i_{r+2}}^{a_{r+2}} \cdots \tau_{i_l}^{a_l} \beta B_{1\,2}^{-1}]\\
&&\quad - [ \tau_2^{a_1} \cdots \tau_2^{a_{r-1}} \tau_1^{a_r}\tau_3^{a_{r+1}} (\sigma_3-(q-1)) 
\tau_{i_{r+2}}^{a_{r+2}} \cdots \tau_{i_l}^{a_l} \beta B_{1\,2}^{-1}]\\
&=&\omega_1 + \omega_2 \in \Span (\DD_{l-1,\infty} \cup \EE_{r+1,l,\infty})\,.
\end{array}\]
\qed

\bigskip\noindent
At this point, thanks to Claims 1 and 2, we have proved that
\begin{equation}\label{R418}
\Span (\DD_{l,\infty}) = \Span (\DD_{l-1,\infty} \cup \EE_{l,l,\infty})\,,
\end{equation}
for all $l \ge 3$.

\bigskip\noindent
{\bf Claim 3.} {\it Let $l \ge 3$. Then
\begin{equation}\label{R419}
\Span( \DD_{l-1,\infty} \cup \EE_{l,l,\infty}) = \Span( \DD_{l-1,\infty} \cup \EE_{l,l,3})\,.
\end{equation}}

\bigskip\noindent
{\bf Proof.} It suffices to show that
\[
\Span( \DD_{l-1,\infty} \cup \EE_{l,l,n}) = \Span( \DD_{l-1,\infty} \cup \EE_{l,l,n-1})
\]
for all $n \ge 4$.

\bigskip\noindent
Set $\varepsilon =2$ if $l$ is even, and $\varepsilon=1$ if $l$ is odd. Let $\omega \in \EE_{l,l,n}$ be 
an element of the form $\omega= [ \tau_1^{a_1} \tau_2^{a_2} \tau_1^{a_3} \cdots \tau_\varepsilon^{a_l} 
\beta]$, where $\beta \in \BB_n$. By construction, either $\beta \in \BB_{n-1}$, or $\beta = \alpha_1 
\sigma_{n-1} \alpha_2$ for some $\alpha_1, \alpha_2 \in B_{n-1}$. If $\beta \in \BB_{n-1}$, then $\omega 
\in \EE_{l,l,n-1}$. If $\beta = \alpha_1 \sigma_{n-1} \alpha_2$ for some $\alpha_1, \alpha_2 \in B_{n-1}$, 
then
\[
\omega= z [ \tau_1^{a_1} \tau_2^{a_2} \tau_1^{a_3} \cdots \tau_\varepsilon^{a_l} \alpha_1 \alpha_2] \in 
\Span (\DD_{l-1,\infty} \cup \EE_{l,l,n-1})\,.
\]
\qed

\bigskip\noindent
{\bf Claim 4.} {\it Let $l \ge 3$. Then
\begin{equation}\label{R420}
\Span (\DD_{l-1,\infty} \cup \EE_{l,l,3}) = \Span( \DD_{l-1,\infty})\,.
\end{equation}}

\bigskip\noindent
{\bf Proof.} Let
\[
\delta_0 = (z^2 -(q-1)z -q)^{-1} (z-(q-1)+\sigma_1)\,.
\]
A direct calculation shows that we have
\[
\delta_0 (z-\sigma_1) = (z-\sigma_1) \delta_0 = 1
\]
in $\K (z) \otimes \HH (B_2)$.

\bigskip\noindent
We set $\varepsilon=1$ if $l$ is odd, and $\varepsilon=2$ if $l$ is even. Let $\omega \in \EE_{l,l,3}$. 
We write $\omega= [\tau_1^{a_1} \tau_2^{a_2} \cdots \tau_\varepsilon^{a_l} \beta]$, where $\beta \in 
\BB_3$. Set
\[\begin{array}{rcl}
\omega_1 &=& [\tau_1^{a_1} \tau_3^{a_2} \tau_2^{a_3} B_{1\,2} \tau_2^{a_4}\cdots \tau_\varepsilon^{a_l} 
\beta \delta_0]\,,\\
\omega_2 &=& [\tau_1^{a_1+a_3} \tau_3^{a_2} (\sigma_3-\sigma_1) \tau_2^{a_4}\cdots 
\tau_\varepsilon^{a_l} \beta \delta_0]\,,\\
\omega_3 &=& [\tau_1^{a_1+a_2} B_{1\,2}\tau_1^{a_3} \tau_2^{a_4}\cdots \tau_\varepsilon^{a_l} \beta 
\delta_0]\,.
\end{array}\]
Obviously, $\omega_2, \omega_3 \in \Span (\DD_{l-1,\infty})$. On the other hand, by Lemma 4.11,
\[\begin{array}{rcl}
&&[(\sigma_3-\sigma_1)^2 \tau_1^{a_1} \tau_3^{a_2} \tau_2^{a_3} B_{1\,2} \tau_2^{a_4}\cdots 
\tau_\varepsilon^{a_l} \beta \delta_0]\\
&=& [(\sigma_3-\sigma_1)^2 \tau_1^{a_1} \tau_2^{a_2+a_3} B_{1\,2} \tau_2^{a_4}\cdots 
\tau_\varepsilon^{a_l} \beta \delta_0]\\
&&\quad + [(\sigma_3-\sigma_1)^2 \tau_3^{a_2} \tau_2^{a_1+a_3} B_{1\,2} \tau_2^{a_4}\cdots 
\tau_\varepsilon^{a_l} \beta \delta_0]\\
&&\quad + [(\sigma_3-\sigma_1) (B_{1\,2} - B_{2\,3}) \tau_2^{a_1+a_2+a_3} B_{1\,2} \tau_2^{a_4}\cdots 
\tau_\varepsilon^{a_l} \beta \delta_0]\\
&\in&\Span (\DD_{l-1,\infty})\,.
\end{array}\]
Moreover, by Lemma 4.10,
\[\begin{array}{rcl}
&&[C_{1\,3} \tau_1^{a_1} \tau_3^{a_2} \tau_2^{a_3} B_{1\,2} \tau_2^{a_4}\cdots \tau_\varepsilon^{a_l} 
\beta \delta_0]\\
&=&[C_{1\,3} \tau_1^{a_1+a_2} \tau_2^{a_3} B_{1\,2} \tau_2^{a_4}\cdots \tau_\varepsilon^{a_l} \beta 
\delta_0]\\
&\in& \Span (\DD_{l-1,\infty})\,.
\end{array}\]
Hence, by Lemma 4.10,
\[\begin{array}{rcl}
\omega_1 &=& (q+1)^{-2} [(\sigma_3-\sigma_1)^2 \tau_1^{a_1} \tau_3^{a_2} \tau_2^{a_3} B_{1\,2} 
\tau_2^{a_4}\cdots \tau_\varepsilon^{a_l} \beta \delta_0]\\
&&\quad + (q+1)^{-2} [C_{1\,3} \tau_1^{a_1} \tau_3^{a_2} \tau_2^{a_3} B_{1\,2} \tau_2^{a_4}\cdots 
\tau_\varepsilon^{a_l} \beta \delta_0]\\
&\in& \Span (\DD_{l-1,\infty})\,.
\end{array}\]
Finally, by \eqref{R42},
\[\begin{array}{rcl}
\omega &=& [ (z- \sigma_1) \tau_1^{a_1} \tau_2^{a_2} \tau_1^{a_3} \tau_2^{a_4} \cdots 
\tau_\varepsilon^{a_l} \beta \delta_0]\\
&=& [ (\sigma_3- \sigma_1) \tau_1^{a_1} \tau_2^{a_2} \tau_1^{a_3} \tau_2^{a_4} \cdots 
\tau_\varepsilon^{a_l} \beta \delta_0]\\
&=& [ \tau_1^{a_1} \sigma_3 \tau_2^{a_2} \tau_1^{a_3} \tau_2^{a_4} \cdots \tau_\varepsilon^{a_l} \beta 
\delta_0] - [  \tau_1^{a_1} \sigma_1 \tau_2^{a_2} \tau_1^{a_3} \tau_2^{a_4} \cdots 
\tau_\varepsilon^{a_l} \beta \delta_0]\\
&=& [ \tau_1^{a_1} \tau_3^{a_2}\sigma_2  \tau_1^{a_3} \tau_2^{a_4} \cdots \tau_\varepsilon^{a_l} \beta 
\delta_0] + [ \tau_1^{a_1+a_3} \tau_3^{a_2}(\sigma_3-(q-1))  \tau_2^{a_4} \cdots \tau_\varepsilon^{a_l} 
\beta \delta_0]\\
&&\quad - [ \tau_1^{a_1+a_2} B_{1\,2} \tau_1^{a_3} \tau_2^{a_4} \cdots \tau_\varepsilon^{a_l} \beta 
\delta_0]\\
&=& \omega_1+ \omega_2 -\omega_3 \in \Span (\DD_{l-1,\infty})\,.
\end{array}\]
\qed

\bigskip\noindent
At this point we have proved that
\[
\Span (\DD_{l,\infty}) = \Span (\DD_{l-1,\infty})
\]
for all $l \ge 3$. This implies that
\begin{equation}\label{R421}
\TR_d = \Span (\DD_{d,\infty}) = \Span (\DD_{2,\infty})\,.
\end{equation}

\bigskip\noindent
Now, let
\begin{multline*}
\FF_1 = \{ [\tau_1^d],[\tau_1^d \sigma_1]\} \\
\cup \{ [\tau_1^a \tau_2^b], [\tau_1^a \tau_2^b \sigma_1],[\tau_1^a \tau_2^b \sigma_2],[\tau_1^a 
\tau_2^b \sigma_1 \sigma_2],[\tau_1^a \tau_2^b \sigma_2 \sigma_1],[\tau_1^a \tau_2^b \sigma_1 \sigma_2 
\sigma_1]\,;\,a,b \ge 1 \text{ and } a+b=d\}\,.
\end{multline*}

\bigskip\noindent
{\bf Claim 5.} {\it $\TR_d = \Span (\FF_1)$.}

\bigskip\noindent
{\bf Proof.} One can easily prove using the same arguments as in the proof of Claim 1 that
\[
\Span (\DD_{2, \infty}) = \Span (\EE_{1,1,\infty} \cup \EE_{2,2,\infty})\,.
\]
On the other hand, using the same arguments as in the proof of Claim 3, it is easily seen that
\[
\Span (\EE_{1,1,\infty} \cup \EE_{2,2,\infty}) = \Span (\FF_1)
\]
\qed

\bigskip\noindent
Now, let
\[
\FF_2= \{ [\tau_1^d], [\tau_1^d \sigma_1]\} \cup \{ [\tau_1^a \tau_2^b], [\tau_1^a \tau_2^b 
\sigma_1]\,;\, a,b \ge 1 \text{ and } a+b=d\}\,.
\]

\bigskip\noindent
{\bf Claim 6.} {\it $\TR_d = \Span (\FF_2)$.}

\bigskip\noindent
{\bf Proof.} Let $a,b \ge 1$ such that $a+b=d$. Then
\begin{multline*}
[\tau_1^a \tau_2^b \sigma_1 \sigma_2 \sigma_1] = [\tau_1^a  \sigma_1 \sigma_2 \sigma_1 \tau_1^b] = 
[\tau_1^{a+b}  \sigma_1 \sigma_2 \sigma_1] = z[\tau_1^{a+b} \sigma_1^2] \\
= z(q-1) [\tau_1^d \sigma_1] + zq[\tau_1^d] \in \Span (\FF_2)\,.
\end{multline*}
\[
[\tau_1^a \tau_2^b \sigma_1 \sigma_2] = [\tau_1^a  \sigma_1 \sigma_2\tau_1^b] = [\tau_1^{a+b}  \sigma_1 
\sigma_2] = z [\tau_1^d  \sigma_1] \in \Span (\FF_2)\,.
\]
By \eqref{R43} we have
\[
[\tau_1^a \tau_2^b B_{1\,2}] = [\tau_1^a B_{1\,2} \tau_1^b] = [\tau_1^{a+b} (\sigma_1+\sigma_2-(q-1))] 
= [\tau_1^d (\sigma_1+z-(q-1))] \in \Span (\FF_2)\,.
\]
On the other hand,
\[
[\tau_1^a \tau_2^b \sigma_2] = [\tau_1^a \tau_2^b B_{1\,2}] - [\tau_1^a \tau_2^b \sigma_1] +(q-
1)[\tau_1^a \tau_2^b]\,,
\]
thus $[\tau_1^a \tau_2^b \sigma_2] \in \Span (\FF_2)$. By \eqref{R44} we have
\[\begin{array}{cl}
&[\tau_1^a \tau_2^b (\sigma_1 \sigma_2 + \sigma_2 \sigma_1 -(q-1) \sigma_1 -(q-1) \sigma_2 + (q^2-
q+1))]\\
=&[\tau_1^{a+b}(\sigma_1 \sigma_2 + \sigma_2 \sigma_1 -(q-1) \sigma_1 -(q-1) \sigma_2 + (q^2-q+1))]\\
=&[\tau_1^d(2z\sigma_1  -(q-1) \sigma_1 -(q-1) z + (q^2-q+1))]\\
\in&\Span (\FF_2)\,.
\end{array}\]
On the other hand,
\begin{multline*}
[\tau_1^a \tau_2^b \sigma_2 \sigma_1] = [\tau_1^a \tau_2^b (\sigma_1 \sigma_2 + \sigma_2 \sigma_1 -(q-
1) \sigma_1 -(q-1) \sigma_2 + (q^2-q+1))]\\
 - [\tau_1^a \tau_2^b \sigma_1 \sigma_2 ] +(q-1) [\tau_1^a \tau_2^b \sigma_1] + (q-1) [\tau_1^a \tau_2^b  
\sigma_2] -(q^2-q+1) [\tau_1^a \tau_2^b]\,,
\end{multline*}
thus $[\tau_1^a \tau_2^b \sigma_2 \sigma_1] \in \Span (\FF_2)$.
\qed

\bigskip\noindent
Let
\begin{multline*}
\FF_3 = \{ [\tau_1^d], [\tau_1^d \sigma_1]\} \cup \{ [\tau_1^a \tau_2^b]\,;\, a \ge b \ge 1 \text{ and 
} a+b=d\} \\
\cup \{ [\tau_1^a \tau_2^b \sigma_1]\,;\, a>b\ge 1 \text{ and } a+b=d\}\,.
\end{multline*}
Note that $|\FF_3| = d+1$, thus the following finishes the proof of Theorem 4.7.

\bigskip\noindent
{\bf Claim 7.} {\it $\TR_d = \Span (\FF_3)$.}

\bigskip\noindent
{\bf Proof.} Let $a,b \ge 1$ such that $a<b$ and $a+b=d$.
\[
[\tau_1^a \tau_2^b] = [\tau_2^b \tau_1^a] = [B_{1\,2} \tau_2^b \tau_1^a B_{1\,2}^{-1}] = [\tau_1^b 
\tau_2^a] \in \Span (\FF_3)\,.
\]
By Lemma 4.8,
\[
[\tau_1^b \tau_2^a B_{1\,2}] = [\tau_1^b  B_{1\,2}\tau_1^a] = [\tau_1^{a+b} (\sigma_1+\sigma_2-(q-1))] 
= [\tau_1^d (\sigma_1+z-(q-1))] \in \Span(\FF_3)\,.
\]
Moreover,
\[
[\tau_1^b \tau_2^a \sigma_2] = [\tau_1^b \tau_2^a B_{1\,2}] - [\tau_1^b \tau_2^a \sigma_1] + (q-
1)[\tau_1^b \tau_2^a]\,,
\]
thus $[\tau_1^b \tau_2^a \sigma_2] \in \Span (\FF_3)$. It follows that
\[
[\tau_1^a \tau_2^b \sigma_1] = [\sigma_1 \tau_1^a \tau_2^b] = [\tau_1^a \sigma_1 \tau_2^b] = [\tau_2^b 
\tau_1^a \sigma_1] = [B_{1\,2} \tau_2^b \tau_1^a \sigma_1 B_{1\,2}^{-1}] = [\tau_1^b \tau_2^a \sigma_2] 
\in \Span (\FF_3)\,.
\]
Now, assume that $d$ is even, and let $a=b=\frac{d}{2}$. We have
\[
[\tau_1^a \tau_2^a \sigma_1] = [\sigma_1 \tau_1^a \tau_2^a] = [\tau_1^a \sigma_1 \tau_2^a] = [\tau_2^a 
\tau_1^a \sigma_1] = [B_{1\,2} \tau_2^a \tau_1^a \sigma_1 B_{1\,2}^{-1}] = [\tau_1^a \tau_2^a 
\sigma_2]\,.
\]
Moreover,
\[
[\tau_1^a \tau_2^a B_{1\,2}] = [\tau_1^a  B_{1\,2}\tau_1^a] = [\tau_1^d (\sigma_1 + \sigma_2 -(q-1))] = 
[\tau_1^d (\sigma_1 + z -(q-1))] \in \Span (\FF_3)\,.
\]
Thus
\[
[\tau_1^a \tau_2^a \sigma_1] = \frac{1}{2}([\tau_1^a \tau_2^a \sigma_1]+[\tau_1^a \tau_2^a \sigma_2]) = 
\frac{1}{2} [\tau_1^a \tau_2^a B_{1\,2}] + \frac{1}{2} (q-1) [\tau_1^a \tau_2^a] \in \Span (\FF_3)\,.
\]
\qed


\section{Universal Markov trace and universal HOMFLY-type invariant}

Let $X,Y$ be two new variables. We define the {\it universal Markov trace} as the collection $\hat T = 
\{ \widehat{\tr}_n\}_{n=1}^{+\infty}$ of $\K$-linear maps
\[
\widehat{\tr}_n : \HH(SB_n) \to \C(\sqrt{q},z)[X,Y]\,, \quad n\ge 1,
\]
defined as follows. Let $d \ge 0$, and let $\omega \in \HH (S_dB_n)$. Then
\[
\widehat{\tr}_n (\omega) = \sum_{k=0}^d \frac{\sqrt{q}^k}{(d-k)!\, k!} X^k\cdot Y^{d-k}\cdot T_k^d (\omega)\,,
\]
where $\{ T_0^d, T_1^d, \dots, T_d^d \}$ is the $\K(z)$-basis of $\TR_d$ constructed in Section 4.

\bigskip\noindent
{\bf Proposition 5.1.}
{\it 
\begin{enumerate}
\item
We have $\widehat{\tr}_n (\alpha \beta) = \widehat{\tr}_n (\beta \alpha)$ for all $\alpha, \beta \in 
SB_n$, and all $n \ge 1$.
\item
Let $\iota_n: \HH(SB_n) \to \HH(SB_{n+1})$ be the morphism induced by the inclusion $SB_n 
\hookrightarrow SB_{n+1}$. Then $\widehat{\tr}_{n+1} \circ \iota_n = \widehat{\tr}_n$ for all $n \ge 
1$.
\item
$\widehat{\tr}_{n+1}(\iota_n(\omega) \sigma_n) = z \cdot \widehat{\tr}_n (\omega)$ for all $\omega \in 
\HH(SB_n)$, and all $n \ge 1$.
\item
We have
\[
\widehat{\tr}_n (\tau_i \omega) = X\sqrt{q} \cdot \widehat{\tr_n} (\sigma_i \omega) + Y
\cdot \widehat{\tr}_n (\omega)
\]
for all $\omega \in \HH (SB_n)$, all $n \ge 2$, and all $1 \le i\le n-1$.
\end{enumerate}}

\bigskip\noindent
{\bf Proof.} Parts (1), (2), and (3) follow from the definition of a Markov trace (see Section 3), and 
from the fact that $T_0^d, T_1^d, \dots, T_d^d$ are Markov traces for all $d \ge 0$.

\bigskip\noindent
We turn now to prove (4). Let $\beta \in S_dB_n$. We write $\beta$ in the form
\[
\beta= \alpha_0 \tau_{i_1} \alpha_1 \cdots \tau_{i_d} \alpha_d\,,
\]
where $1 \le i_j \le n-1$ for $1 \le j\le d$, and $\alpha_j \in B_n$ for $0 \le j \le d$. For $S 
\subset \{1, \dots, d\}$ we set
\[
\beta (S)= \alpha_0 u_1 \alpha_1 \cdots u_d \alpha_d\,,
\]
where $u_j= \sigma_{i_j}$ if $j \in S$, and $u_j=1$ if $j \not\in S$. It is easily checked that, for $0 
\le k \le d$, $T_k^d (\beta)$ is given by the formula
\begin{equation}\label{R51}
T_k^d (\beta)= k!\, (d-k)! \sum_{\substack{S \subset \{1, \dots, d\}\\ |S|=k}} T_0^0 (\beta(S))\,.
\end{equation}
This implies that
\begin{equation}\label{R52}
\widehat{\tr}_n (\beta)= \sum_{S \subset \{1, \dots, d\}}\sqrt{q}^{|S|} X^{|S|} Y^{d-|S|} \cdot T_0^0 (\beta(S))\,.
\end{equation}
Now, from \eqref{R52}, it follows that
\begin{align*}
\widehat{\tr}_n (\tau_i \beta) &= X \sqrt{q} \sum_{S \subset \{1, \dots, d\}} \sqrt{q}^{|S|}X^{|S|} 
Y^{d-|S|} \cdot T_0^0 (\sigma_i 
\beta(S)) + Y \sum_{S \subset \{1, \dots, d\}} \sqrt{q}^{|S|} X^{|S|} Y^{d-|S|}\cdot T_0^0 (\beta(S))\\
&= X\sqrt{q} \cdot \widehat{\tr}_n (\sigma_i \beta) + Y \cdot\widehat{\tr}_n (\beta)\,.
\end{align*}
\qed

\bigskip\noindent
Recall from Section 3 that $\pi: SB_n \to \HH (SB_n)$ denotes the natural map, and that $\varepsilon: 
SB_n \to \Z$ is the homoùmorphism defined by
\[
\varepsilon(\sigma_i)=1\,, \quad \varepsilon(\sigma_i^{-1}) =-1\,, \quad \varepsilon(\tau_i) =0\,, 
\quad \text{for } 1 \le i\le n-1\,.
\]
We consider the following change of variables:
\[
z= \frac{q-1}{1-qy} \quad \Leftrightarrow \quad y= \frac{z-q+1}{qz}\,.
\]
For $\beta \in SB_n$, we set
\[
\hat I (\beta) = \left( \frac{q-1}{1-qy} \right)^{-n+1} \cdot (\sqrt{y})^{\varepsilon (\beta) -n+1} 
\cdot \widehat{\tr}_n (\pi (\beta))\,.
\]
This is an element of $\C(\sqrt{q},\sqrt{y}) [X,Y]$.

\bigskip\noindent
The following can be proved in the same way as Proposition 3.3.

\bigskip\noindent
{\bf Proposition 5.2.} {\it Let $(\alpha,n)$ and $(\beta,m)$ be two singular braids. If $\hat \alpha$ 
is isotopic to $\hat \beta$, then $\hat I (\alpha) = \hat I (\beta)$.}
\qed

\bigskip\noindent
Let $\LL$ denote the set of (isotopy classes of) singular links. For $L \in \LL$, we choose a singular 
braid $(\beta,n)$ such that $\hat \beta =L$, and we set $\hat I (L)= \hat I (\beta)$. By Proposition 
5.2, the map $\hat I: \LL \to \C(\sqrt{q},\sqrt{y}) [X,Y]$ is a well-defined invariant that we call the {\it 
universal HOMFLY-type invariant} of $\LL$.

\bigskip\noindent
For $d \ge 0$, we denote by $\SS_d$ the set of invariants $I: \LL_d \to \C (\sqrt{q}, \sqrt{y})$ which 
satisfies the skein relation for $t= \sqrt{y} \sqrt{q}$ and $x= \sqrt{q} - \frac{1}{\sqrt{q}}$. Now, 
the above terminology ``universal HOMFLY-type invariant'' is justified by the following.

\bigskip\noindent
{\bf Theorem 5.3.} {\it Let $d \ge 0$, and let $L,L' \in \LL_d$. We have $\hat I (L) = \hat I(L')$ if 
and only if $I(L)=I(L')$ for all $I \in \SS_d$.}

\bigskip\noindent
{\bf Proof.} Let $\TR_d'$ be the space of traces on $\{ \HH (S_dB_n)\}_{n=1}^{+\infty}$ with 
coefficients in $\C (\sqrt{q}, \sqrt{y})$. Clearly, $\TR_d'$ is a $\C(\sqrt{q}, \sqrt{y})$-vector 
space, and
\[
\TR_d'= \C(\sqrt{q}, \sqrt{y}) \otimes \TR_d\,.
\]
On the other hand, we have
\begin{equation}\label{R53}
\hat I (\beta) = \sum_{k=0}^d \frac{\sqrt{q}^k}{(d-k)!\, k!} X^k Y^{d-k}\cdot I_{T_k^d} (\beta)\,,
\end{equation}
for all $\beta \in S_dB_n$.

\bigskip\noindent
Let $L,L' \in \LL_d$ such that $\hat I (L)= \hat I (L')$. By \eqref{R53}, we have $I_{T_k^d} (L) = 
I_{T_k^d} (L')$ for all $0 \le k\le d$. Let $I \in \SS_d$. By Proposition 3.5, there exists $T \in 
\TR_d'$ such that $I=I_T$. By Theorem 4.7, there exist $\lambda_0, \lambda_1, \dots, \lambda_d \in 
\C(\sqrt{q}, \sqrt{y})$ such that
\[
T= \lambda_0 T_0^d + \lambda_1 T_1^d + \cdots + \lambda_d T_d^d\,.
\]
Then
\[
I(L)= \sum_{k=0}^d \lambda_k I_{T_k^d} (L) = \sum_{k=0}^d \lambda_k I_{T_k^d} (L') = I(L')\,.
\]

\bigskip\noindent
Now, let $L,L' \in \LL_d$ such that $I(L)=I(L')$ for all $I \in \SS_d$. We have in particular 
$I_{T_k^d} (L)= I_{T_k^d} (L')$ for all $0 \le k\le d$, thus, by \eqref{R53}, $\hat I(L) = \hat I 
(L')$.
\qed

\bigskip\noindent
Let $A$ be an abelian group, let $I: \LL \to A$ be an invariant, and let $X,Y \in A$. We say that $I$ 
satisfies the {\it $(X,Y)$ desingularization relation} if
\[
I(L_X)= X \cdot I(L_+) + Y \cdot I(L_0)\,,
\]
for all singular links $L_X, L_+, L_0 \in \LL$ that have the same link diagram except in the 
neighborhood of a crossing where they are like in Figure 5.1.

\begin{figure}[htb]
\bigskip
\centerline{
\setlength{\unitlength}{.5cm}
\begin{picture}(15,4)
\put(0,1){\includegraphics[width=7.5cm]{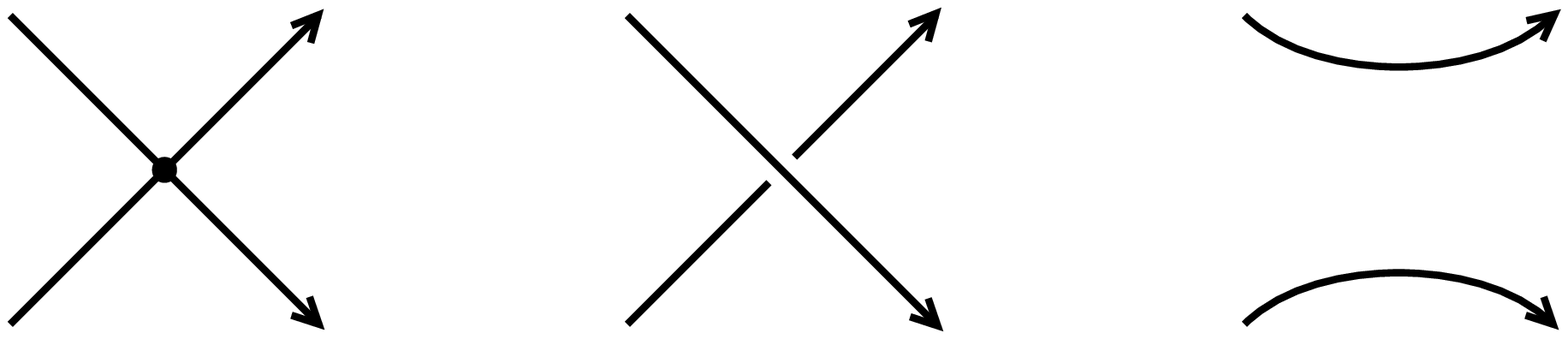}}
\put(1,0){\small $L_X$}
\put(7,0){\small $L_+$}
\put(13,0){\small $L_0$}
\end{picture}}
\bigskip
\centerline{{\bf Figure 5.1.}  The singular links $L_X$, $L_+$, and $L_0$.}
\end{figure}

\bigskip\noindent
We set
\[
t= \sqrt{y} \sqrt{q}\,, \quad x= \sqrt{q} - \frac{1}{\sqrt{q}}\,,
\]
and we define $\widetilde{\tr}_n : SB_n \to \C (\sqrt{q}, \sqrt{y}) [X,Y]$ by
\[
\widetilde{\tr}_n (\beta)= (\sqrt{y})^{-\varepsilon (\beta)} \cdot \widehat{\tr}_n (\pi (\beta))\,.
\]
With these new notations, $\hat I (\beta)$ can be written
\[
\hat I (\beta) = \left( \frac{1-t^2}{tx}\right)^{n-1} \cdot t^{\varepsilon (\beta)} \cdot \widetilde{\tr}_n 
(\beta)\,.
\]

\bigskip\noindent
{\bf Proposition 5.4.} {\it The invariant $\hat I$ satisfies the $(t,x)$ skein relation and the $(X, Y)$ 
desingularization relation.}

\bigskip\noindent
{\bf Proof.} The fact that $\hat I$ satisfies the $(t,x)$ skein relation is proved in the same way as 
Proposition~3.4. So, we only need to show that $\hat I$ satisfies the $(X, Y)$ desingularization 
relation.

\bigskip\noindent
Let $L_X, L_+, L_0 \in \LL$ be three singular links that have the same link diagram except in the 
neighborhood of a crossing where they are like in Figure 5.1. A careful reading of the proof of 
Theorem 2.3 shows that there exist a singular braid $(\beta,n)$ and an index $1 \le i\le n-1$ such that 
$L_X = \widehat{ \tau_i \beta}$, $L_+ = \widehat{\sigma_i \beta}$, and $L_0= \hat \beta$. On the other 
hand, Proposition 5.1.(4) implies that
\[
\widetilde{\tr}_n (\tau_i \beta) = Xt \cdot \widetilde{\tr}_n (\sigma_i \beta) + Y\cdot \widetilde{\tr}_n 
(\beta)\,.
\]
Hence
\[\begin{array}{rcl}
X \cdot \hat I (L_+) + Y\cdot \hat I (L_0) &=& \left( \frac{1-t^2}{tx} \right)^{n-1} \cdot t^{\varepsilon( 
\tau_i \beta)} \cdot (Xt \cdot \widetilde{\tr}_n (\sigma_i \beta) +Y \cdot \widetilde{\tr}_n (\beta))\\
&=& \left( \frac{1-t^2}{tx} \right)^{n-1} \cdot t^{\varepsilon( \tau_i \beta)} \cdot \widetilde{\tr}_n 
(\tau_i \beta)\\
&=& \hat I(L_X)\,.
\end{array}\]
\qed

\bigskip\noindent
Now, the following shows that our invariant $\hat I$ is a reasonable extension of the HOMFLY 
polynomial to the singular links.

\bigskip\noindent
{\bf Theorem 5.5.} {\it There exists a unique invariant $\hat I : \LL \to \C (\sqrt{q}, \sqrt{y}) 
[X,Y]$ which satisfies the $(t,x)$ skein relation and the $(X,Y)$ desingularization relation, and 
which takes the value $1$ on the trivial knot. Moreover, $\hat I (L) \in \C [t^{\pm 1}, x^{\pm 1}, 
X,Y]$ for all $L \in \LL$.}

\bigskip\noindent
{\bf Proof.} The existence of the invariant is given by Proposition 5.4.

\bigskip\noindent
Suppose that $\hat I' : \LL \to \C (\sqrt{q}, \sqrt{y}) [X,Y]$ is an invariant which satisfies the 
$(t,x)$ skein relation and the $(X,Y)$ desingularization relation, and which takes the value 1 on 
the trivial knot. Let $L \in \LL_d$ be a singular link with $d$ singular points. We prove by 
induction on $d \ge 0$ that $\hat I' (L)= \hat I(L)$, and that this element belongs to $\C [t^{\pm 1}, 
x^{\pm 1}, X,Y]$.

\bigskip\noindent
The case $d=0$ is well-known (see \cite{Jones1}, \cite{HOMFLY}). We assume $d \ge 1$. Let $P$ 
be a singular point of $L$. Set $L_X=L$, and let $L_+$ and $L_0$ be the singular links having the same 
link diagram as $L$ except in the neighborhood of $P$ where they are like in Figure 5.1. Then, by 
induction and by the $(X,Y)$ desingularization relation, we have
\[
\hat I' (L) = X \cdot \hat I' (L_+) +Y \cdot \hat I' (L_0) = X \cdot \hat I(L_+) +Y \cdot \hat I(L_0) = \hat I(L)\,.
\]
Moreover, again by induction,
\[
\hat I(L) = X \cdot \hat I (L_+) +Y \cdot \hat I(L_0) \in \C [t^{\pm 1}, x^{\pm 1}, X,Y]\,.
\]
\qed



\bigskip\bigskip\noindent
{\bf Luis Paris},

\smallskip\noindent
Institut de Math\'ematiques de Bourgogne, UMR 5584 du CNRS, Universit\'e de Bourgogne, B.P. 
47870, 21078 Dijon cedex, France

\smallskip\noindent
E-mail: {\tt lparis@u-bourgogne.fr} 

\bigskip\noindent
{\bf Lo\"{\i}c Rabenda},

\smallskip\noindent 
Institut de Math\'ematiques de Bourgogne, UMR 5584 du CNRS, Universit\'e de Bourgogne, B.P. 
47870, 21078 Dijon cedex, France

\smallskip\noindent
E-mail: {\tt lrabenda@u-bourgogne.fr}


\begin{thebibliography}{99}

\bibitem{Baez1}
{\bf J.C. Baez.}
{\it Link invariants of finite type and perturbation theory.}
Lett. Math. Phys. {\bf 26} (1992), no. 1, 43--51.

\bibitem{Birma1}
{\bf J.S. Birman.}
{\it New points of view in knot theory.}
Bull. Amer. Math. Soc. (N.S.) {\bf 28} (1993), no. 2, 253--287.

\bibitem{HOMFLY}
{\bf P. Freyd, D. Yetter, J. Hoste, W.B.R. Lickorish, K. Millett, A. Ocneanu.}
{\it A new polynomial invariant of knots and links.}
Bull. Amer. Math. Soc. (N.S.) {\bf 12} (1985), no. 2, 239--246.

\bibitem{Garsi1}
{\bf F.A. Garside.}
{\it The braid group and other groups.}
Quart. J. Math. Oxford Ser. (2) {\bf 20} (1969), 235--254.

\bibitem{Gemei1}
{\bf B. Gemein.}
{\it Singular braids and Markov's theorem.}
J. Knot Theory Ramifications {\bf 6} (1997), no. 4, 441--454.

\bibitem{Jones2}
{\bf V.F.R. Jones.}
{\it A polynomial invariant for knots via von Neumann algebras.}
Bull. Amer. Math. Soc. (N.S.) {\bf 12} (1985), no. 1, 103--111.

\bibitem{Jones1}
{\bf V.F.R. Jones.}
{\it Hecke algebra representations of braid groups and link polynomials.}
Ann. of Math. (2) {\bf 126} (1987), no. 2, 335--388.

\bibitem{Kauff1}
{\bf L.H. Kauffman.}
{\it Invariants of graphs in three-space.}
Trans. Amer. Math. Soc. {\bf 311} (1989), no. 2, 697--710.

\bibitem{KauVog1}
{\bf L.H. Kauffman, P. Vogel.}
{\it Link polynomials and a graphical calculus.}
J. Knot Theory Ramifications {\bf 1} (1992), no. 1, 59--104.

\bibitem{PrzTra1}
{\bf J.H. Przytycki, P. Traczyk.}
{\it Invariants of links of Conway type.}
Kobe J. Math. {\bf 4} (1988), no. 2, 115--139.

\end{thebibliography}
\end{document}